\documentclass[12pt]{amsart}
\usepackage{amsmath}\usepackage{amsrefs}\usepackage{wasysym}\usepackage{marvosym}
\usepackage{amsfonts}\usepackage{mathrsfs}\usepackage{verbatim}
\usepackage{amssymb}
\usepackage{amstext}\usepackage{stmaryrd} \usepackage{yfonts}
\usepackage{amsbsy}

\usepackage{amsopn}
\usepackage{amsthm}\usepackage{color}\newtheorem{theorem}{Theorem}
\DeclareMathAlphabet{\mathpzc}{OT1}{pzc}{m}{it}

\newtheorem{example}{Example}
\newtheorem{exple}{Example}
\newtheorem{lem}{Lemma}
\newtheorem{propn}{Proposition}
\newtheorem{rem}{Remark}
\numberwithin{equation}{section}
\numberwithin{theorem}{section}
\numberwithin{propn}{section}
\numberwithin{lem}{section}

\def\undersetbrace#1\to#2{\underbrace{#2}_{#1}}

\def\text#1{\mbox{#1}}
\def\flushpar{\par\noindent}

\def\mod{\mbox{ mod }}

\newcommand{\mapright}[1]{%
    \smash{\mathop{%
        \hbox to 1cm{\rightarrowfill}
        }
    \limits^{#1}
    }
}
\newcommand{\mapleft}[1]{%
    \smash{\mathop{%
        \hbox to 1cm{\rightarrowfill}
        }
    \limits_{#1}
    }
}

\def\e{\epsilon}
\def\a{\alpha}
\def\b{\beta}

\def\g{\gamma}
\def\d{\delta}
\def\D{\Delta}
\def\s{\sigma}
\def\Si{\Sigma}
\def\th{\theta}
\def\l{\lambda}
\def\x{\times}

\def\o{\overline}
\def\f{\flushpar}

\def\v{\varphi}

\def\om{\omega}
\def\Om{\Omega}
\def\B{\mathcal B}
\def\T{\widehat T}
\def\({\biggl(}
\def\){\biggr)}
\def\<{\langle}
\def\>{\rangle}
\def\bdy{\partial}
\def\bul{\smallskip\f$\bullet\ \ \ $}

\def\){\biggr)}

\def\pde{\text{pointwise dual ergodic}}\def\st{\text{such that}}

\def\<{\bold\langle}
\def\>{\bold\rangle}
\def\pprime{\prime\prime}

\def\bul{\smallskip\f$\bullet\ \ \ $}
\def\sms{\smallskip\f}

\def\lra{\longrightarrow}

\def\pf{\smallskip\f{\it Proof} }

\begin{document}

\title{Limit theory for some positive\\ stationary processes with infinite mean}
\author{ Jon. Aaronson,\   Roland Zweim\"{u}ller}
 \address[Jon. Aaronson]{\ \ School of Math. Sciences, Tel Aviv University,
69978 Tel Aviv, Israel.}

\email{aaro@tau.ac.il}
\address[Roland Zweim\"{u}ller]{\ \ Faculty of Mathematics, University of Vienna,
1090 Vienna, Austria }
\email{rzweimue@member.ams.org}
\begin{abstract} We prove distributional limit theorems and one-sided laws of
the iterated logarithm for a class of positive, mixing, stationary, stochastic processes
which contains those obtained from non-integrable observables over certain piecewise
expanding maps. This is done by extending Darling-Kac theory to a suitable family of
infinite measure preserving transformations.
\end{abstract}

\maketitle\markboth{J. Aaronson \& R. Zweim\"{u}ller}{Limit theory for some positive stationary processes}

\subsection*{\sl Overview}
We prove  limit theorems for positive, stationary, processes (SPs) with infinite mean satisfying
mixing conditions which occur naturally in certain dynamical systems:
 {\it Stable limit theorems}  (SLTs)   for certain $\vartheta_\mu$-mixing SPs and
 {\it One-sided laws of the iterated logarithm} (LILs) for certain $\psi^*$-mixing SPs (definitions below).

 The method of proof is by  {\it inversion} which is done
 by first  building a {\it Kakutani tower} over the generating probability
 preserving transformation, using the time zero observation as height function.

The mixing properties of the stationary process ensure that the resulting infinite
measure preserving transformation is {\it weakly pointwise dual ergodic},
which allows us to develop a generalized {\it Darling Kac theory} for ergodic sums of
this system. The results for the original stochastic process then follow by
a standard inversion argument.

We illustrate both the finite-measure and the infinite-measure results
by applying them to certain one-dimensional dynamical systems.

\section{Definitions and Preliminaries}

\subsection*{\sl Stationary processes}
We are going to consider partial sums of ergodic $\Bbb R_+$-valued stationary
processes $(\xi_n)_{n\ge0}$ with $\Bbb E(\xi_n)=\infty$. Such a process can always be
represented as $\xi_n=\v\circ S^n$, where $S$ is a measure preserving
transformation ({\it m.p.t.}) on a probability space $(\Om,\mathcal A ,P)$,
and $\v:\Om\to\Bbb R_+$ is measurable with $\Bbb E(\v)=\infty$
(and w.l.o.g. $\mathcal A=\sigma\{\v\circ S^n : n\ge 0\}$).
Due to non-integrability, it will suffice to restrict attention to the $\Bbb N$-valued
case, as by the ergodic theorem the partial sums of the fractional parts are
asymptotically negligible compared to the partial sums of the process.
If indeed $\v:\Om\to\Bbb N$, we let $\a=\a_\v=\{[\v=l]:\ l\in \Bbb N\}$,
and thus obtain a {\it probability preserving fibred system}
$(\Om,\mathcal A ,P,S,\a)$ in the sense of the following definition.

\subsection*{\sl Nonsingular transformations and Fibred systems}
A measurable map $S$ on a $\sigma$-finite space $(\Om,\mathcal A ,m)$
is called {\it nonsingular} if $m \circ S^{-1} \ll m$. Its {\it  transfer operator}
(w.r.t. $m$) is the positive linear map $\widehat S:L^1(m)\to L^1(m)$ defined by
$$\int_A\widehat S fdm=\int_{S^{-1}A} fdm\ \ \ \ (f\in L^1(m),\ A\in\mathcal A).$$
A {\it fibred} or {\it piecewise (pcw) invertible system} is a quintuple
 $(\Om,\mathcal A ,m,S,\a)$ where $S$ is a nonsingular transformation on
 $(\Om,\mathcal A ,m)$, and $\a\subset\mathcal A $ is a countable, unilateral
 generator so that the restriction $S:a\to Sa$ is invertible, non-singular
 on each $a\in\a$. In this case, for every $k\ge 1$, $(\Om,\mathcal A ,m,S^k,\a_k)$ is a
fibred system, where $\a_k:=\bigvee_{j=0}^{k-1}S^{-j}\a$.

The transfer operator of $(\Om,\mathcal A ,m,S,\a)$ can be represented as
$$\widehat S f=\sum_{a\in\a}1_{Sa}v_a' (f\circ v_a),$$
where $v_a:Sa\to a$ denotes the inverse of $S:a\to Sa$, and $v_a':=\tfrac{dm\circ v_a}{dm}$.

If $m$ actually is an $S$-invariant probability measure,
the system is called {\it probability preserving}, and we write $P:=m$.

\subsubsection*{\sl Mixing}
We let $\mathcal P(\Om,\mathcal A )$ denote the collection of
probability measures on $(\Om,\mathcal A )$, and call $\mu\in\mathcal P(\Om,\mathcal A )$
\emph{equivalent} to $P$, $\mu\sim P$ if $\mu \ll P \ll \mu$.
The probability preserving fibred system  $(\Om,\mathcal A ,P,S,\a)$ is called\ \bul {\it
$\vartheta_\mu$-mixing}\  (for some $\mu\sim P$) if $\vartheta_\mu(n)\to 0$, where
\begin{align*}\vartheta_\mu(n) :=\sup\,\left\{\tfrac{|P(A\cap S^{-(n+k)}B)-P(A)P(B)|}{\mu(B)}:\ k\ge 1,
A\in\s(\a_k),\ B\in\mathcal A\right\};\end{align*}
\bul {\it reverse $\phi$-mixing} if $\phi_-(n)\to 0$, where
\begin{align*}\phi_-(n) :=\sup\,\left\{\tfrac{|P(A\cap S^{-(n+k)}B)-P(A)P(B)|}{P(B)}:\ k\ge 1,
A\in\s(\a_k),\ B\in\mathcal A\right\};\end{align*}
\bul {\it $\psi^*$-mixing}
if $\psi^*(n)\to 1$, where
\begin{align*}\psi^*(n) :=\sup\,\left\{\tfrac{P(A\cap S^{-(n+k)}B)}{P(A)P(B)}:\ k\ge 1,
A\in\s(\a_k),\ B\in\mathcal A\right\};\end{align*}

\bul {\it $\psi$-mixing} if $\psi(n)\to 0$, where
\begin{align*}\psi(n) :=\sup\,\left\{\tfrac{|P(A\cap S^{-(n+k)}B)-P(A)P(B)|}{P(A)P(B)}:\ k\ge 1,
A\in\s(\a_k),\ B\in\mathcal A\right\};\end{align*}
and {\it continued fraction mixing} if, in addition to $\psi$-mixing, $\psi(1)<\infty$.

\begin{rem}
\

{\bf a)}  As shown in [{Br1}], $\psi^*(1)<\infty$  implies $\psi^*$-mixing.
Elementary computation shows that  $\phi_-(n)\le\psi^*(n)-1$ so  $\psi^*$-mixing entails reverse $\phi$-mixing.
Note that  $\psi^*(1)\le 1+\psi(1)$. For  examples with $\psi^*(1)<\infty$ which are not $\psi$-mixing,\ see chapter 5 in [{Br2}].

{\bf b)} Note that $\vartheta_P\equiv\phi_-$. In \S6, we consider  a class of interval maps (weakly mixing {\tt RU}   maps)  for which
$\vartheta_\mu(n)\to 0$ exponentially.  For these interval maps (as shown in [AN]) $\psi^*(1)<\infty$ implies continued fraction mixing (see \S6).
\end{rem}

\subsubsection*{\sl Strong distributional convergence and Limit laws}
For $(X,\B,m)$  a $\s$-finite measure space, $F_n:X\to [0,\infty]$ measurable, and $Y \ge0$ a random variable,
we say that {\it $(F_n)$ converges strongly in distribution to $Y$}, written
$$F_n\overset{\mathfrak d}{\underset{n\to\infty}\lra } Y ,$$
if it converges in law with respect to all absolutely continuous probabilities, that is, if
$$\int_Xg(F_n)dP\underset{n\to\infty}\lra \, \Bbb E(g(Y))\
\ \forall\ \ g\in \mathcal C ([0,\infty]),\ P\in\mathcal P(X,\B),\ P\ll m.$$

\

For  $\g\in [0,1]$ we let $Y_\g\ge0$ denote a random variable which has the
{\it normalized Mittag-Leffler distribution of order $\g$}, that is,
$\Bbb E(Y_\g^p)=\tfrac{p!(\Gamma(1+\gamma))^p}{\Gamma(1+p\gamma)}$ for $p\ge0$.
Evidently $Y_1\equiv 1$, and $Y_0$ has exponential distribution.
Also, $Y_{\frac12}$ is the absolute value of a centered Gaussian random variable.

For $\g\in (0,1]$, the variable $Z_\g:=Y^{-\frac1\g}$ then has a
{\it positive $\g$-stable distribution} with
$\Bbb E(e^{-tZ_\g})=\exp({-\Gamma(1+\gamma)\,t^\g)}$ for $t>0$.

\section{Results on stationary processes}

 In the statements below, $(\Om,\mathcal A,P,S,\a)$ is a
 probability preserving fibred system, and
 $\v:\Om\to\Bbb N$ is $\a$-measurable. We let
$$ \v_n:=\sum_{k=1}^n\v\circ S^k, \;\;\;\;  n\ge 0,$$
denote the partial sums of the stationary process
$(\xi_n)_{n\ge0}=(\v\circ S^n)_{n\ge0}$, and define
$$ a(n):=\sum_{k=1}^nP([\v_k\le n]),  \;\;\;\;  n\ge 0.$$
In order to establish our results, we'll need to assume that the growth  of $a(n)$ is adapted  to the decay of the mixing coefficients of the process. The main condition is as follows although we need a stronger version (2.6) in Theorem \ref{T_osLIL} (below).
\subsection*{Definition of Adaptedness} \ \ Let $\tau(n)\downarrow\ 0$. We'll say that the increasing sequence $(a(n))_{n\in\Bbb N}$ is {\it adapted to $(\tau(n))_{n\in\Bbb N}$} if
\begin{equation}\label{eq_uaaah}
\frac{n\tau(\delta a(n))}{a(n)}\underset{n\to\infty}\lra\ 0
\;\;\;\text{ for all } \delta>0.
\end{equation}
Note that if  $(a(n))_{n\in\Bbb N}$ is  regularly varying with positive index, then it is adapted to $(\tau(n))_{n\in\Bbb N}$   as soon as (2.1) holds for one $\delta>0$.

\

Our first result is a distributional limit theorem.
In the barely infinite measure case (${\g}=1$) it
comes with an associated a.e. result. For
${\g} \in (0,1)$, corresponding statements will be
established under stronger assumptions in Theorem
\ref{T_osLIL} below.

\begin{theorem}\label{T_SLT}
Suppose that $(\Om,\mathcal A,P,S,\a)$ is a  $\vartheta_\mu$-mixing
probability preserving fibred system, and that $\v:\Om\to\Bbb N$ is $\a$-measurable.
Let $a(n)$ be $\g$-regularly varying with $\gamma \in (0,1]$ and adapted to $(\vartheta_\mu(n))_{n\in\Bbb N}$.
\sms {\bf(a) (Stable limit theorem)}

\begin{equation}
\frac{\v_n}{b(n)}\ \overset{\mathfrak d}{\underset{n\to\infty}\lra} \ Z_\g,
\end{equation}
where $b$ is asymptotically inverse to $a$ in that
$b(a(n)) \sim a(b(n)) \sim n$
(and hence $\tfrac{1}{\g}$-regularly varying).
\sms {\bf(b) (One-sided law of the iterated logarithm for ${\g}=1$)}
\

If, in addition, ${\g}=1$ and $b(n/\log\log n) \log\log n \sim b(n)$,
then
\begin{equation}\label{eq_LIL1}
\varliminf_{n\to\infty}\frac{\v_n}{b(n)}= 1\ \
 \text{a.s.}
\end{equation}

\end{theorem}

\begin{rem}
\

{\bf a)} Theorem 2.1(a) was established for  $\phi$-mixing processes in [{S}]
and for continued fraction mixing processes in [{D}] (see also [{A3}]).

{\bf b)} The functional version of (a) is also valid, and can be proved using
a straightforward, appropriate adaptation of [{B}].

{\bf c)} Theorem 2.1(b) was established for  $\psi$-mixing processes in [{AD1}].

\end{rem}

The results mentioned in  Remark a) also compute the $a(n)$
from the marginal distributions, for which additional
``{\it close correlation}" assumptions such as
$\psi^*(1)<\infty$ are required.
We now show how to determine the asymptotics
of $a(n)$ from the marginal distributions under the weaker
close correlation condition (\ref{eq_Wheelchair}) (but we
still use the stronger $\psi^*(1)<\infty$ in Theorem \ref{T_osLIL} below).

\

\begin{theorem}[{\bf Identifying the normalization}]\label{T_NormalizationID}
Let $(\Om,\mathcal A,P,S,\a)$ be a  $\vartheta_\mu$-mixing
probability preserving fibred system, that $\v:\Om\to\Bbb N$ is $\a$-measurable,
and that there exists some $\Psi\in L^1(P)_+$ such that
\begin{equation}\label{eq_Wheelchair}
\widehat S(\v\wedge n)\le \Psi \, \Bbb E(\v\wedge n)\ \forall\ n\ge 1.
\end{equation}
Assume that $A$ is strictly increasing, regularly varying with index
$\gamma\in (0,1]$, and adapted to $(\vartheta_\mu(n))_{n\in\Bbb N}$,
then
\begin{equation}\label{eq_NormIDcondition3}
\Bbb E(\v\wedge n) \underset{n\to\infty}\sim\ \tfrac{n}{\Gamma(2-\gamma)\Gamma(1+\gamma) A(n)}
\end{equation}
implies
$$a(n) \underset{n\to\infty}\sim A(n).$$
\end{theorem}

Finally, under stronger assumptions, we establish the
following pointwise result, where $C_\gamma:= K_\gamma^{-1/\gamma}$
with
$K_\gamma :=\tfrac {\Gamma(1+\gamma )}{\gamma ^\gamma (1-\g )^{1-\gamma }}$
for $\g\in (0,1)$.

\begin{theorem}[{\bf The one-sided law of the iterated logarithm}]\label{T_osLIL}
\sms   Suppose that $(\Om,\mathcal A,P,S,\a)$ is a $\psi^*$-mixing
probability preserving fibred system with $\psi^*(1)<\infty$,
and that   $\v:\Om\to\Bbb N$ is $\a$-measurable.
\

If $a(n)$ is $\g$-regularly varying for some $\g\in (0,1)$, and
\begin{equation}
\frac{n\phi_-(\delta a(a(n)))}{a(n)}\underset{n\to\infty}\lra 0\
\;\;\;\text{ for some } \delta>0,
\end{equation}
then, for any sequence $(\tau(n))$ with
$\tau(n)\uparrow$ and $\tau(n)/n\downarrow$ as
$n \to \infty$,
\begin{equation} \varliminf_{n\to\infty}\frac{\v_n}{b(n/\tau(n))\tau(n)}
\ge C_\gamma\ \
 \text{a.s.\ \;\;\;\; if} \ \sum_{n=1}^\infty\tfrac 1ne^{-\beta\tau(n)}
<\infty\ \; \forall \beta>1,\end{equation}

\begin{equation} \varliminf_{n\to\infty}\frac{\v_n}{b(n/\tau(n))\tau(n)}
\le C_\gamma\ \
 \text{a.s.\ \;\;\;\; if} \ \sum_{n=1}^\infty\tfrac 1ne^{-r\tau(n)}
=\infty\ \; \forall r<1,\end{equation}
and
\begin{equation} \varliminf_{n\to\infty}\frac{\v_n}{b(n/\log\log(n))
\log\log(n)}
= C_\gamma\ \
 \text{a.s.}\end{equation}
\end{theorem}

\begin{rem}
The conclusion of Theorem \ref{T_osLIL} was established for  iid processes in [{W}],
and for $\psi$-mixing processes in [{AD2}]. Our proof of Theorem \ref{T_osLIL} is by
establishing the conditions needed for the methods of [{AD2}].
Therefore the functional version also follows as in [{AD3}].
\end{rem}

\subsubsection*{\sl Inversion: Kakutani towers and Return time processes}
\

Our results  will be established using the well-known technique of ``inverting" corresponding results for infinite
measure preserving transformations, the connection being established via the
following concept.
The {\it Kakutani tower} of $(\Om,\mathcal A ,P,S,\v)$ is the  object $(X,\B,m,T)$ with $(X,\B,m)$ the
$\s$-finite space defined by
\bul $X:=\bigcup_{n\ge 1}[\v\ge n]\x\{n\}$,
\bul $\B:=\{\bigcup_{n\ge 1}B_n\x\{n\}:\ B_n\in\mathcal A\cap [\v\ge n]\ \forall\ n\ge 1\}$,
\bul $m(A\x\{n\}):=P(A)$,

\noindent and $T:X\to X$ is the map given by
\bul $T(x,n):=\begin{cases} & (x,n+1)\ \ \ \ \ \v(x)>n,\\ & (Sx,1)\ \ \ \ \ \ \ \ \ \v(x)=n.\end{cases}$

\noindent It follows that $(X,\B,m,T)$ is a conservative, measure preserving
transformation which is ergodic iff $(\Om,\mathcal A ,P,S)$ is ergodic.\\

This ``tower building process" is reversible. Given a conservative ergodic
measure preserving system $(X,\B,m,T)$ we define the
{\it return time  process of  $T$ on
$\Om\in\mathcal F:=\{B\in\B: 0<m(B)<\infty\}$} as the
$\Bbb N$-valued stationary process $(\v_\Om\circ T_\Om^n)_{n\ge0}$
on $(\Om,\mathcal B\cap\Om ,m_\Om,)$, where
\bul $\v_\Om(x):=\min\,\{n\ge 1:\ T^nx\in\Om\}$,
\bul $T_\Om(x):=T^{\v(x)}(x)$, and
\bul $m_\Om(A):=m(A\cap\Om)/m(\Om).$

\noindent It follows that the Kakutani tower of $(\Om,\mathcal B\cap\Om ,m_\Om,T_\Om,\v_\Om)$
is a factor of $(X,\B,m',T)$ where $m'=\tfrac1{m(\Om)}m$ (and an isomorph in case $T$ is invertible).\\

Now set $\v_j:=\sum_{i=0}^{j-1}\v_\Om\circ T_\Om^i$, which is the time
of the $n$th return to $\Om$. It is straightforward to check that these
are dual to the {\it occupation times} of $\Om$,
$S_n(1_\Om):= \sum_{k=0}^{n-1} 1_\Om \circ T^k$ in that
\begin{equation}\label{eq_duality}
S_n(1_\Om) \le j \;\;\; \text{ iff } \;\;\; \v_j \ge n.
\end{equation}
This entails, via routine arguments, that various properties of $(\v_j)_{j\ge1}$
are equivalent to corresponding properties of $(S_n(1_\Om))_{n\ge1}$.
Specifically, suppose that $a(n)$ is $\g$-regularly varying with $\g\in (0,1]$,
and let $b$ be (asymptotically) inverse to $a$. Then, for $Y \ge 0$  a random variable,
\begin{equation}\label{eq_duality1}
\tfrac1{a(n)}S_n(1_\Om)\overset{\mathfrak d}\lra m(\Om)\,Y \;\;\; \text{ iff } \;\;\;
\tfrac{\v_{n}}{b(n)}\overset{\mathfrak d}\lra (\tfrac1{m(\Om)Y})^{\frac1\g},
\end{equation}
and
\begin{equation}\label{eq_duality2}
\underset{n\to\infty}\varlimsup\tfrac1{a(n)}S_n(1_\Om)\overset{\text{\tiny a.e.}}=m(\Om) \;\;\; \text{ iff } \;\;\; \underset{n\to\infty}\varliminf\tfrac{\v_{n}}{b(n)}\overset{\text{\tiny a.e.}}=(\tfrac1{m(\Om)})^{\frac1\g}.
\end{equation}

\section{Weak pointwise dual ergodic measure preserving transformations}

In this section, we consider the properties of infinite ergodic systems needed in the proofs of the results of the previous section.

\subsection*{\sl Weak, pointwise dual ergodicity}
Let $T$ be a conservative, ergodic, measure
preserving transformation (not necessarily invertible) on
the $\sigma$-finite space $(X,\B,m)$,
and $\T:L^1(m)\to L^1(m)$ its transfer operator, which
naturally extends to all non-negative measurable functions.
Invariance of $m$ means that $\T 1_X=1_X$, and since $T$ is c.e., any
measurable $g:X\to [0,\infty)$ which is subinvariant, $\T g\le g$,
is, in fact, constant.
Hurewicz's ratio ergodic theorem (Theorem 2.2.1 of [{A1}]), guarantees that
\begin{equation}\label{eq_Hu}
\frac{\sum_{k=0}^{n-1} \T^k f}{\sum_{k=0}^{n-1} \T^k g}
\underset{n\to\infty}{\lra}
\frac{m(f)}{m(g)}\ \;\;\;\; \text{a.e. on}\ X
\end{equation}
for all $f,g \in L^1_+(m):=\{f\in L^1(m): f\ge0\ \text{and}\ m(f)>0 \}$.
(Due to conservativity, $\sum_{k=0}^{n-1} \T^k f \to \infty $ a.e. for such $f$.)\\

Throughout, convergence in measure, ${\overset{m}{\lra}}$,
for our $\sigma$-finite measure $m$, is understood to mean convergence
in measure, ${\overset{\nu}{\lra}}$, for every finite $\nu \ll m$
(or, equivalently, for all $\nu=m_A$ with $A\in \mathcal{F}$).\\

The c.e.m.p.t. $(X,\B,m,T)$ will be called {\it weakly $\pde$}\ if
there exist constants $ a_n=a_n(T)>0, n\ge 1,$ such that
\begin{equation}\label{eq_wp_i}
\frac1{a_n}\sum_{k=0}^{n-1}\T^k f \underset{n\to\infty}{\overset{m}{\lra}}
\int_Xfdm\ \;\;\;\; \text{for}\ f\in L^1_+(m),
\end{equation}
and
\begin{equation}\label{eq_wp_ii}
\varlimsup_{n\to\infty}\frac1{a_n}\sum_{k=0}^{n-1}\T^k f =\int_Xfdm\ \;\;\;\;
\text{a.e. on}\ X\ \;\; \text{for}\ f\in L^1_+(m).
\end{equation}
This generalizes the notion of {\it pointwise dual ergodicity}
(cf. \S 3.7 of [{A1}], or [{A2}]),
which requires $a_n=a_n(T)>0$ such that
\begin{equation}\label{eq_pde}
\frac1{a_n}\sum_{k=0}^{n-1}\T^k f \underset{n\to\infty}{\lra} \int_Xfdm\
\;\;\;\; \text{a.e. on}\ X\ \;\; \text{for}\ f\in L^1_+(m).
\end{equation}

\begin{rem}
No \emph{invertible} c.e.m.p.t. $(X,\B,m,T)$ with $m(X)=\infty$ is pointwise
dual ergodic. (Since in this case $\T f = f \circ T^{-1}$, so that (\ref{eq_pde})
would give a pointwise ergodic theorem with normalizing constants $a_n$ for
$T^{-1}$, which is impossible, see \S 2.4 of [{A1}].) However,
invertible systems can still be weakly pointwise dual ergodic.

For a concrete example, let $T:(0,1)\to(0,1)$ be given by $Tx:=\frac{x}{1-x}$ for
$x<\frac{1}{2}$ and $Tx:=2x-1$ for $x>\frac{1}{2}$, which
is conservative ergodic w.r.t. the invariant measure $m$ with density
$\frac{1}{x}$, and define $a_n:=n/\log n$. By the Darling-Kac theorem for pointwise
dual ergodic transformations, (see [DK], [{A2}], \S 3.7 of  [{A1}],
or [{Z3}]),
$$ \frac1{a_n}\sum_{k=0}^{n-1} f \circ T^k \underset{n\to\infty}{\overset{m}{\lra}}
\int_Xfdm\ \;\;\;\; \text{for}\ f\in L^1_+(m),$$
and according to Proposition 2 of [AD1],
$$\varlimsup_{n\to\infty}\frac1{a_n}\sum_{k=0}^{n-1} f \circ T^k =\int_Xfdm\ \;\;\;\;
\text{a.e. on}\ X\ \;\; \text{for}\ f\in L^1_+(m).$$
It is easily seen that these carry over to the natural extension $T_*$ of $T$.
Therefore, the invertible c.e.m.p.t. $T_*^{-1}$, with transfer operator
$f \mapsto f \circ T_*$ is weakly pointwise dual ergodic
(and hence so is $T_*$).
\end{rem}

Conditions of this flavour can be exploited most efficiently if one succeeds
in identifying special sets on which there is additional control on
the convergence. Recall (cf. [{A1}], [T2]) that $A\in\mathcal F=\{B\in\B: 0<m(B)<\infty\}$ is called
a {\it uniform set} ({\it for} $f\in L^1_+(m)$), written $A\in \mathcal{U}(T)$, if there are
$ a_n=a_n(T)>0$, s.t. (mod $m$)
$$\sup_A \left| \frac1{a_n}\sum_{k=0}^{n-1}\T^k f-m(f) \right|
\underset{n\to\infty}\lra 0,$$
and, more specifically, a {\it Darling-Kac set}, denoted $A\in \mathcal{DK}(T)$, if
$$\sup_A \left| \frac1{a_n}\sum_{k=1}^n\T^k1_A-m(A) \right| \underset{n\to\infty}\lra 0.$$\\
We now define $A\in\mathcal F$ to be a  {\it limited set}, $A\in \mathcal{L}(T)$, if there exist constants
$a_n=a_n(T)>0$ such that
\begin{equation}\label{eq_ls_i}
\frac1{a_n}\sum_{k=0}^{n-1}\T^k 1_A  \underset{n\to\infty}{\overset{m_A}{\lra}} m(A),
\end{equation}
and (mod $m$)
\begin{equation}\label{eq_ls_ii}
\sup_A\, \frac1{a_n}\sum_{k=0}^{n-1}\T^k1_A   \underset{n\to\infty}\lra  m(A).
\end{equation}

In each case the {\it return sequence} $(a_n(T))_{n\ge1}$, which is determined up to asymptotic
equivalence and satisfies $a_{n+1}(T) \sim a_n(T) \to \infty$, can (and will) be taken non-decreasing.
Due to (\ref{eq_Hu}) it does not depend on $f$ or $A$.
It is easy to see that defining
$$  a_n(A):= \sum_{k=0}^{n-1} \frac{m(A\cap T^{-k}A)}{m(A)^2}
\;\;\;\; \text{for}\ A\in \mathcal{F}, $$
we have
\begin{equation}\label{eq_RS_A}
a_n(T) \underset{n\to\infty}\sim a_n(A)
\;\;\;\; \text{for}\ A\in \mathcal{L}(T).
\end{equation}
The existence of uniform sets is equivalent to pointwise dual ergodicity
(but we do not know if $\mathcal{DK}(T)\ne\emptyset$ for every $\pde$ transformation).
Similarly, weak pointwise dual ergodicity is equivalent to the existence of limited sets,
as we have

\begin{propn}[{\bf Limited sets and weak pde from local behaviour}]\label{Prop_WPDE_local}
Let $T$ be a conservative ergodic m.p.t. on $(X,\mathcal B,m)$.\medskip

\noindent{\bf a)}\ Suppose there are $A\in \mathcal{F}$, $f\in L^1_+(m)$, and
constants $ a_n>0$, $n\ge 1$, such that
\begin{equation}\label{eq_wp_iA}
\frac1{a_n}\sum_{k=0}^{n-1}\T^k f\underset{n\to\infty}{\overset{m_A}{\lra}}\int_Xfdm,
\end{equation}
and
\begin{equation}\label{eq_wp_iiA}
\varlimsup_{n\to\infty}\, \frac1{a_n}\sum_{k=0}^{n-1}\T^k f =\int_Xfdm\ \;\;\;\;
\text{\rm a.e. on}\ A.
\end{equation}
Then, for every $\eta>0$, $T$ possesses a limited set
$A' \in \mathcal{F}_A$ with return sequence $(a_n)_{n\ge 1}$ and $m(A')\ge m(A)-\eta$.\medskip

\noindent{\bf b)}\ If $T$ has a limited set, then it is weakly pointwise dual ergodic.
\end{propn}

\pf \;
{\bf a)}\ Let $B_0:=A$.
By Hurewicz's ratio ergodic theorem we may assume w.l.o.g. that $f=1_{B_0}$,
and given any set $B_j \in \mathcal{F} \cap A$ with $m(B_j)>m(A)-\eta$, we also have
$a_n^{-1} \sum_{k=0}^{n-1} \T^k 1_{B_j} \to m(B_j)$ a.e. on $A$.
Egorov's theorem then provides us with some $B_{j+1} \in \mathcal{F} \cap {B_j}$
such that $m(B_{j+1}) > m(A)-\eta$, and
$$ \frac{1}{a_n} \sum_{k=0}^{n-1} \T^k 1_{B_j}
\underset{n\to\infty}{\lra} m(B_j)\ \;\;\;\text{uniformly on}\ B_{j+1}.$$
Using this to inductively define a decreasing sequence $(B_j)_{j \ge 0}$ in $\mathcal{F}$, we
obtain a set $A' := \bigcap_{j \ge 0} B_j$ with $m(A') \ge m(A)-\eta$.
Given $\e>0$ choose $j$ with $m(B_j)<m(A')+\e/2$, then
$$ \sup_{A'} \frac{1}{a_n} \sum_{k=0}^{n-1} \T^k 1_{A'}
\le \sup_{B_{j+1}} \frac{1}{a_n} \sum_{k=0}^{n-1} \T^k 1_{B_j}
< m(B_j)+\frac{\e}{2}\ \;\;\;\text{for}\ n \ge n_j(\e).$$
This gives the required control from above. Since, by (\ref{eq_Hu}), (\ref{eq_wp_iA})
also holds with $f=1_{A'}$, we see that $A'$ is indeed a limited set.\\

\noindent{\bf b)}\ Now start from the assumption that $T$ has a limited set $A$.
Due to Hurewicz's theorem, weak pde follows as soon as we check
the defining conditions (\ref{eq_wp_i}) and (\ref{eq_wp_ii}) for $f=1_A$, i.e. we have to prove
\begin{equation}\label{eq_ipp}
\frac1{a_n}\sum_{k=0}^{n-1}\T^k 1_A \underset{n\to\infty}{\overset{m}{\lra}} m(A),
\end{equation}
and
\begin{equation}\label{eq_iipp}
\bar{f} := \varlimsup_{n\to\infty}\, \frac1{a_n}\sum_{k=0}^{n-1}\T^k 1_A
= m(A)\ \;\;\;\; \text{\rm a.e. on}\ X.
\end{equation}\\

Letting $A_n:=A\cap[\v_A>n]$, $n\ge 0$, we have, by routine arguments,
\begin{equation}\label{eq_B1}
\sum_{n\ge0} \T^n 1_{A_n} = 1_X,
\end{equation}
and decomposing, for any $n\ge 0$, $A$ according to the time of the last return
before time $n$, $A= A_n \cup \bigcup_{k=0}^{n-1} A\cap T^{-(n-k)}(A_k)$ (disjoint),
we find that for $N\ge0$,
\begin{equation}\label{eq_B2}
\sum_{n=0}^N \T^n 1_{A} =
\sum_{k=0}^{N-1} \T^k \left( 1_{A_k} \sum_{j=0}^{N-k} \T^j 1_A \right)
+ \sum_{n=0}^N \T^n 1_{A_n}.
\end{equation}
Since $A\in \mathcal{L}(T)$, there is some $M\in (0,\infty)$ such that
\begin{equation}\label{eq_B3}
\frac{1}{a_n} \sum_{k=0}^{n-1} \T^k 1_A \le M\ \;\;\;\; \text{on}\ A\ \text{for}\ n\ge 1.
\end{equation}
With (\ref{eq_B3}) and (\ref{eq_B1}) providing bounds for the sums in (\ref{eq_B2}), we get
\begin{equation}\label{eq_B4}
\frac{1}{a_N} \sum_{n=0}^{N-1} \T^n 1_A \le M'\cdot 1_X\ \;\;\;\; \text{for}\ N\ge 1,
\end{equation}
where $M':=M+1/a(1)\in (0,\infty)$. Consequently, $\bar{f} \le M'\cdot 1_X$.
By a canonical version of Fatou's lemma for positive operators, we have
$\T \bar{f} \ge \bar{f}$. Hence $g:=M'\cdot 1_X - \bar{f} \ge 0$ is subinvariant,
$\T g \le g$. As $T$ is conservative ergodic, this implies that $g$ is constant a.e.,
and hence so is $\bar{f}$. In view of (\ref{eq_ls_i}), this yields (\ref{eq_iipp}).\\

To finally prove convergence in measure on $X$, it suffices to check it on
each of the sets $T^{-l}A$, $l\ge0$, since these cover $X$. But for each $l$,
\begin{equation*}
\int_{T^{-l}A} \left( \frac{1}{a_n} \sum_{k=0}^{n-1} \T^k 1_A \right) dm
=  \frac{a_{n+l}(A)-a_{l}(A)}{a_n}\, m(A) \underset{n\to\infty}{\lra} m(A).
\end{equation*}
Together with (\ref{eq_B4}) and (\ref{eq_iipp}) this entails
$a_n^{-1} \sum_{k=0}^{n-1} \T^k 1_A
{\overset{m_{T^{-l}A}}{\lra}} m(A)$,
and hence (\ref{eq_ipp}).\; \qed\\

The notion of return sequences has originally been introduced in the
context of rationally ergodic transformations (cf. \S 3.3 of [{A1}]).
The present use of the term is justified by

\begin{propn}[{\bf Weak pde implies rational ergodicity}]
Let $T$ be a weakly pointwise dual ergodic c.e.m.p.t. on $(X,\mathcal B,m)$.
Then every limited set $A$ satisfies a R\'enyi inequality, meaning that
there is some $M=M(A)\in (0,\infty)$ such that
$$\int_A (S_n(1_A))^2 dm \le M \left(  \int_A S_n(1_A) dm \right)^2\ \;\;\;\;
\text{for}\ n\ge 1.$$
In particular, $T$ is rationally ergodic with $(a_n(T))_{n\ge1}$ a
return sequence in the sense of [{A1}].
\end{propn}

\pf \; Same as the proof of Proposition 3.7.1 in [{A1}], using (\ref{eq_RS_A}) and
the existence of limited sets established above.\ \qed\\

\subsubsection*{\sl (Weak) pointwise dual ergodicity and special sets for Kakutani towers}
Kakutani towers above $\vartheta_\mu$-mixing systems satisfying our adaptedness  conditions are weakly pointwise dual ergodic. In the
presence of regular variation, $\vartheta_\mu$-mixing with sufficiently
fast rate  implies pointwise dual ergodicity.

Recall that if $(X,\B,m,T)$ is a c.e.m.p.t. and $\Om\in\B,\ m(\Om)=1$, has
return time $\v=\v_\Om$, and $\v_j:=\sum_{i=0}^{j-1} \v \circ T_\Om^i$, then,
for $n\ge1$,
\begin{equation}\label{eq_AnThroughPhi}
a_n(\Om)=\sum_{k=1}^n m(\Om\cap T^{-k}\Om)= \sum_{j=1}^n m(\Om\cap [\v_j \le n]).
\end{equation}

\
\begin{theorem}[{\bf (Weak) pde via $\vartheta_\mu$-mixing return processes}]\label{Thm_WPDE_Thetamu}\
Let $(X,\B,m,T)$ be a c.e.m.p.t. and suppose that $\Om\in\B,\ m(\Om)=1$ has a countable partition $\a \subset \mathcal B\cap \Om$ $\st$  $\v=\v_\Om$ is  $\a$-measurable and that $(\Om,\mathcal B\cap \Om,m,T_\Om,\a)$  is $\vartheta_\mu$-mixing for some $\mu \sim m_\Om$.\medskip
\sms {\bf(a)} If the $a_n:=a_n(\Om)$ is adapted to $(\vartheta_\mu(n))_{n\in\Bbb N}$,
then $(X,\B,m,T)$ is weakly $\pde$ with $a_n(T) \sim a_n$,
and \sms $\forall\ \e>0\ \exists\ N_\e\ \ \st\ \forall\  n \ge N_\e,$
\begin{equation}\label{eq_Radioactivity}
\frac1{a_n}\sum_{k=1}^n\T^k 1_\Om -1 \le \e(1+\tfrac{d\mu}{dP})\ \ \text{\rm a.e. on }\ \Om.
\end{equation}
In particular, if $\|\tfrac{d\mu}{dP}\|_\infty<\infty$, then $\Om\in\mathcal L(T)$.\medskip
\sms {\bf(b)}    If $(a_n)_{n\ge1}$ is $\g$-regularly varying for some $\g\in (0,1]$ and
$$\vartheta_\mu(n)=O(\tfrac1{n^r})\ \;\;\;\; \text{for some}\ r>\tfrac1\g-1,$$
then $a_n:=a_n(\Om)$ is adapted to $(\vartheta_\mu(n))_{n\in\Bbb N}$, and $(X,\B,m,T)$ is
pointwise dual ergodic.\medskip
\sms {\bf(c)}    If  $(a_n)_{n\ge1}$ is $\g$-regularly varying for some $\g\in (0,1]$,
$\psi^*(1)<\infty$, and
\begin{equation}\label{eq_CompareAnTheta2}
\frac{n\vartheta_\mu(\delta a(a(n)))}{a(n)}\underset{n\to\infty}\lra 0\ \ \ \text{for all}\  \delta>0,
\end{equation}
then $a_n:=a_n(\Om)$ is adapted to $(\vartheta_\mu(n))_{n\in\Bbb N}$ is satisfied, and $(X,\B,m,T)$ is
pointwise dual ergodic.
Moreover, $\forall\ \e>0\ \exists\ N_\e\ \ \st\ \forall\ n \ge N_\e$
\begin{equation}\label{eq_Radioactivity2}
\left| \frac1{a_n}\sum_{k=1}^n\T^k 1_\Om -1 \right| \le \e(1+\tfrac{d\mu}{dP})\ \ \text{\rm a.e. on }\ \Om.
\end{equation}
In particular, if $\|\tfrac{d\mu}{dP}\|_\infty<\infty$, then $\Om\in\mathcal{DK}(T)$.
\end{theorem}

\pf\ \
We write $(\Om,\mathcal A ,P,S,\a):=(\Om,\mathcal B\cap \Om,m,T_\Om,\a)$.
Note first that $\vartheta_\mu$-mixing implies, for $j,p \ge 1$ and $A \in \alpha_j$,
\begin{equation}
\widehat{S}^{j+p} 1_A   \le  P(A) + \vartheta_\mu(p) \, \tfrac{d\mu}{dP}\
\ \;\;\;\;\; \text{\rm a.e. on}\ \Om,
\end{equation}
as is easily seen by integrating over arbitrary $B \in \mathcal{B} \cap \Om$.
We claim that for $n, p\ge 1$,
\begin{equation}\label{eq_pointer}
\T_n := \sum_{k=1}^n\T^k 1_\Om \le p+a_n+n \vartheta_\mu(p) \, \tfrac{d\mu}{dP}\ \;\;\;\;\; \text{\rm a.e. on}\ \Om.\end{equation}
To see this, observe that $\T^k 1_\Om = \sum_{j \ge 1} \widehat{S}^j 1_{[\v_j = k]}$ a.e. on $\Om$,
and hence
\begin{align*} \T_n &= \sum_{j=1}^n\widehat{S}^j 1_{[\v_j\le n]}
\le\sum_{j=1}^{n+p}\widehat{S}^j 1_{[\v_j\le n]}\\
&\le p+\sum_{j=1}^n\widehat{S}^{j+p} 1_{[\v_{j+p}\le n]}\le p+\sum_{j=1}^n\widehat{S}^{j+p} 1_{[\v_{j}\le n]}\\
&\le  p+\sum_{j=1}^n(P([\v_{j}\le n])+\vartheta_\mu(p)\tfrac{d\mu}{dP})\\  
&=p+a_n+n\vartheta_\mu(p)\tfrac{d\mu}{dP}\ \;\;\;\;\; \text{a.e. on}\ \Om,
\end{align*}
since $[\v_{j}\le n]\in\a_j$.\\

\pf {\it of } (a):  To establish (\ref{eq_Radioactivity}), we let
$p = p_{n,\e}:=[\e a_n]$ for $\e>0$. By adaptedness we can choose $N_\e$ $\st$
$n\vartheta_\mu(p_{n,\e})<\e a_n$ whenever $n \ge N_\e$. Then by (\ref{eq_pointer}),
\begin{equation}\label{eq_goodN}
\T_n \le a_n(1+\e(1+\tfrac{d\mu}{dP}))\ \;\;\;\;\; \text{a.e. on}\ \Om\ \text{for}\ n\ge N_\e,
\end{equation}
as required.\\

To prove weak pointwise dual ergodicity, we will use
part a) of Proposition \ref{Prop_WPDE_local}.
Write $R_n:= \T_n / a_n$
and observe first that by (\ref{eq_Radioactivity}),
\begin{equation}\label{eq_lhd}\
\varlimsup_{n\to\infty} R_n \le 1\ \;\;\;\;\; \text{a.e. on}\ \Om.
\end{equation}
In fact, (\ref{eq_Radioactivity}) shows that for any $M>0$, the estimate (\ref{eq_lhd})
holds uniformly on $\Om \cap [\tfrac{d\mu}{dP}\le M]$. Consequently,
\begin{equation}\label{eq_triangleright}\
P (\left[ R_n \ge 1+\e' \, \right]) \underset{n\to\infty}{\lra} 0\ \;\;\;\;\; \text{for every}\ \e'>0.
\end{equation}

To obtain the analogous statement for convergence from below,
we observe that for $t \in (0,1)$,
\begin{equation}\label{eq_UsedTwice}
R_n \le 1+\e(1+\tfrac{d\mu}{dP})\; \text{a.e. on } \Om\ \; \Longrightarrow \;
P([R_n\le t]) \le\tfrac{2\e}{1+\e-t},
\end{equation}
since (${\mathbb E}$ denoting expectation w.r.t. $P$)
\begin{align*}
1 &={\mathbb E}(R_n)={\mathbb E}(R_n1_{[R_n>t]})+{\mathbb E}(R_n1_{[R_n\le t]})\\
  &\le (1+\e)P([R_n>t])+\e\mu([R_n>t])+t P([R_n\le t])\\
  &\le 1+2\e-(1+\e-t)P([R_n\le t]).
\end{align*}
Now fix $\e\in (0,1)$ and take $N_\e$ as in (\ref{eq_goodN}).
For $t:=1-\sqrt\e \in (0,1)$, observation (\ref{eq_UsedTwice}) then yields
\begin{equation*}
P([R_n\le 1-\sqrt\e \, ])\le 2\sqrt\e\ \;\;\;\;\; \text{for}\  n\ge N_\e.
\end{equation*}
This readily implies
\begin{equation*}\
P (\left[ R_n \le 1-\e' \, \right]) \underset{n\to\infty}{\lra} 0\ \;\;\;\;\; \text{for every}\ \e'>0,
\end{equation*}
which, together with (\ref{eq_triangleright}), gives (\ref{eq_wp_iA}).
Combined with (\ref{eq_lhd}) the latter yields (\ref{eq_wp_iiA}). (There is some subsequence
$n_k \nearrow \infty$ such that $R_{n_k} \to 1$ a.e. for $P$.)\\

\pf {\it of } (b):
It is easily seen that the present assumptions entail adaptedness.
As in part (a) this implies (\ref{eq_lhd}).

To check pointwise dual ergodicity, it remains to show that
\begin{equation}\label{eq_rhd}
\varliminf_{n\to\infty}R_n \ge 1 \ \;\;\;\;\; \text{a.e. on}\ \Om.
\end{equation}
(Proposition 3.7.5 of [A1] ensures that a.e. convergence on some $A\in \mathcal{F}$ suffices).
To this end, choose
$c \ge 1$ such that $\vartheta_\mu(n) \le c / n^r$ for all $n \ge 1$.
There exist $s>0$ and $N_0$ so that
$$a(n)>n^{\frac1{r+1}+2s}\ \ \forall\ n>N_0.$$
Choosing $p=n^{\frac1{r+1}}$ in (\ref{eq_pointer}) we have
$$R_n=\tfrac{\T_n}{a(n)}\le 1+\tfrac1{n^{2s}}(1+\tfrac{d\mu}{dP})\ \text{on}\ \Om\ \forall\ n\ge N_0.$$
Due to (\ref{eq_UsedTwice}), we then see that for $t \in (0,1)$,
$$q_t(n) := P([R_n \le t]) \le \tfrac{2c}{n^{s}} \cdot (1-t+\tfrac{c}{n^{s}})^{-1}
  \le \tfrac{2c}{(1-t) n^{s}}\ \;\;\;\;\; \text{for}\ n \ge N_0.$$
Since, for all $t$ and $\l>1$, $\sum_{n\ge 1} q_t([\l^n])<\infty$, BCL now implies
\begin{equation}\label{eq_rhdrhd}
\varliminf_{n\to\infty}R_{[\l^n]}\ \ge\ 1\ \;\;\;\;\; \text{a.s.\ on}\ \Om\ \text{for all}\ \l>1.
\end{equation}
To finally prove convergence (\ref{eq_rhd}) of the full sequence,
fix any $\l>1$ and choose integers $\kappa_n(\l)\nearrow \infty$ so that
$[\l^{\kappa_n(\l)}]\le n\le [\l^{\kappa_{n}(\l)+1}]$. Then regular variation of $(a_n)_{n\ge1}$ yields
$$\frac{\T_n}{a_n}\ge \frac{\T_{[\l^{\kappa_{n}(\l)}]}}{a_{[\l^{\kappa_{n}(\l)+1}]}}\sim \frac1{\l^\g} \, \frac{\T_{[\l^{\kappa_{n}(\l)}]}}{a_{[\l^{\kappa_{n}(\l)}]}}\ \;\;\;\;\; \text{a.e. on}\ \Om\ \text{as}\ n\to\infty.$$
In view of (\ref{eq_rhdrhd}), we thus have
$\varliminf_{n\to\infty}R_{n}\ \ge\ \l^{-\g}$ a.s.\ on $\Om$ for all $\l>1$, and (\ref{eq_rhd}) follows.\\

\pf {\it of } (c): Note first that adaptedness, and hence (\ref{eq_Radioactivity}) holds.
Thus,  to prove (\ref{eq_Radioactivity2}), it suffices to check that for $n \ge N_\e,$
\begin{equation}\label{eq_Biohazard}
\tfrac1{a(n)}\sum_{k=1}^n\T^k 1_\Om\ge 1-\e(1+\tfrac{d\mu}{dP})\ \ \; \text{a.e. on } \Om.
\end{equation}
We show first that $\forall\ p,q,n\in\Bbb N$ with $n \ge q$,
\begin{equation}\label{eq_Pointinghand}
\T_n\ge a(n)-p-n\vartheta_\mu(p)\tfrac{d\mu}{dP}-\psi^*(1)^2(P([\v_{p}>q])a(n)+q)\
\end{equation}
a.e. on  $\Om$. To see this, observe that
\begin{align*}
\T_n&  =\sum_{k=1}^n\widehat{S}^j 1_{[\v_j\le n]}\\ &\ge \sum_{j=1}^n\widehat{S}^{j+p} 1_{[\v_{j+p}\le n]}-p\\ &= \sum_{j=1}^n\widehat{S}^{j+p} 1_{[\v_{j}\le n]}-\sum_{j=1}^n\widehat{S}^{j+p} 1_{[\v_j\le n<\v_{j+p}]}-p\\ &=:\Si_1-\Si_2-p
\end{align*}
Now, because $[\v_{j}\le n]\in\a_j$, we have
\begin{align*}
\Si_1 \ge \sum_{j=1}^n(P([\v_{j}\le n])-\vartheta_\mu(p)\tfrac{d\mu}{dP})\ \ =a(n)-n\vartheta_\mu(p)\tfrac{d\mu}{dP}.
\end{align*}
On the other hand,
\begin{align*}
\Si_2 &\le
\psi^*(1)\sum_{j=1}^nP([\v_j\le n<\v_{j+p}])\\
&=\psi^*(1)\sum_{j=1}^n\sum_{\ell=j}^nP([\v_j=\ell,\ \v_{p}\circ S^j>n-\ell])\\ &\le
\psi^*(1)^2\sum_{j=1}^n\sum_{\ell=j}^nP([\v_j=\ell])P([\v_{p}>n-\ell])\\
&=:\psi^*(1)^2(\Si_2'+\Si_2'')
\end{align*}
with
\begin{align*}
\Si_2' &:=\sum_{j=1}^n\sum_{\ell=j}^{n-q}P([\v_j=\ell])P([\v_{p}>n-\ell])\\ &
\le P([\v_{p}>q])\sum_{j=1}^n\sum_{\ell=j}^{n-q}P([\v_j=\ell])\\ &\le
P([\v_{p}>q])\sum_{j=1}^nP([\v_j\le n])\\ &= P([\v_{p}>q])a(n);
\end{align*}
and
\begin{align*}
\Si_2''&:=\sum_{j=1}^n \sum_{\ell=n-q+1}^nP([\v_j=\ell])P([\v_{p}>n-\ell])\\
&\le \sum_{j=1}^n \sum_{\ell=n-q+1}^nP([\v_j=\ell])\\
&= \sum_{j=1}^n P([n-q<\v_j\le n])\\
&= \sum_{k=n-q+1}^nm(\Om\cap T^{-k}\Om)\le q.
\end{align*}
Putting this together establishes (\ref{eq_Pointinghand}).\\

To finally check (\ref{eq_Biohazard}) for some given $\e>0$, choose
$$ \d \in (0, \tfrac{\e}{3\psi^*(1)^2})\;\;  \text{for which}\;\; \Pr(Z_\g>\tfrac1\d)<\tfrac{\e}{3\psi^*(1)^2},$$
and let
$$ p_n=p_{n,\e}:=[\d^{2\g} a(a(n))],\;\; \text{and}\;\;
   q_n=q_{n,\e}:=[\d a(n)].$$
Then $a^{-1}(p_n)\sim \d^2a(n)$, and $\tfrac{q_n}{a^{-1}(p_n)}\underset{n\to\infty}\lra\tfrac1\d$.
Under the present assumptions, the stable limit theorem, Theorem \ref{T_SLT},
applies to our $(\v_n)$, as its proof below does not depend on part c) of Theorem \ref{Thm_WPDE_Thetamu}. Therefore,
$$P([\v_{p_n}>q_n])\underset{n\to\infty}\lra \Pr(Z_\g>\tfrac1\d)<\tfrac{\e}{3\psi^*(1)^2}.$$
Now choose $N_\e$  so large that, for all $n\ge N_\e$,
$$p_n<\tfrac{\e}{3}a(n), \;\; \ n\vartheta_\mu(p_{n})<\e a(n)\ \;\; \text{and}\;\; \ P([\v_{p_n}>q_n])<\tfrac{\e}{3\psi^*(1)^2}. $$
Then, using (\ref{eq_Pointinghand}), we find indeed that for $n\ge N_\e$,
\begin{align*}
\T_n &\ge a(n)-p_n-n\vartheta_\mu(p_n)\tfrac{d\mu}{dP}-\psi^*(1)^2(P([\v_{p_n}>q_n])a(n)+q_n)\\ &\ge
a(n)-\tfrac{\e}{3}a(n)-\e a(n)\tfrac{d\mu}{dP}-\tfrac{2\e}3a(n)\\ &=(1-\e(1+\tfrac{d\mu}{dP}))a(n)
\end{align*}
a.e. on $\Om$, which is (\ref{eq_Biohazard}).
By (\ref{eq_Radioactivity}) and (\ref{eq_Biohazard}),
$$\frac1{a(n)}\sum_{k=1}^n\T^k1_{\Om}\underset{n\to\infty}\lra 1 \ \ \; \text{a.e. on } \Om,$$
so that $T$ is pointwise dual ergodic (Proposition 3.7.5 of [{A1}] again). This convergence
is in fact uniform on each $\Om \cap [\tfrac{d\mu}{dP}\le M]$, $M>0$,
whence the assertion about Darling-Kac sets. \;\;$\Box$

\section{Moment sets and the stable limit theorem}

\subsubsection*{\sl Darling-Kac theorem and stable limits}
The statement of the stable limit theorem announced above is dual to a
Darling-Kac type limit theorem for the Kakutani tower, which we now
establish in the setup of weakly pointwise dual ergodic systems.
It quantifies ``{return rates}" and
determines the limit distribution of occupation
times $S_n(1_A)=\sum_{k=0}^{n-1}1_{A}\circ T^k$ of sets $A$ of finite measure:\medskip

\begin{theorem}[{\bf Darling-Kac theorem for weakly pde systems}]\label{T_DK}
Let $T$ be a weakly  pointwise dual ergodic c.e.m.p.t. on $(X,\B,m)$.
If its return sequence $(a_n(T))_{n\ge 1}$ is regularly varying of
index $\a \in [0,1]$, then
$$\frac{S_n(f)}{a_n(T)}  \overset{\mathfrak d}{\underset{n\to\infty}\lra}\  m(f)\,Y_\g\ \;\;\;\; \text{for}\ f\in L_+^1(m),$$
where  $Y_\g$ has the normalised Mittag-Leffler distribution of
order $\g$.
\end{theorem}

The stable limit theorem follows easily from this:

\subsection*{Proof of Theorem \ref{T_SLT}(a)}
Let $(X,\B,m,T)$ be the {\it Kakutani tower} of $(\Om,\mathcal A,P,S,\v)$, then  $(X,\B,m,T)$ is weakly $\pde$ with return sequence $a(n)$. By the Darling-Kac Theorem,
$$\tfrac1{a(n)}S_n(1_\Om)\ \overset{\mathfrak d}\lra\  Y_\g,$$
and we need only recall (\ref{eq_duality1}). \qed

\subsubsection*{\sl Moment sets}

The proof of the Darling-Kac theorem identifies
sets for which the asymptotics of the moments of the occupation time distributions
can be understood. Let $(X,\B,m,T)$ be a conservative,
ergodic, measure preserving transformation.
For $A\in\mathcal F=\{F\in\B,\ 0<m(F)<\infty\}$, recall that
$$ a_n(A)=\sum_{k=0}^{n-1}\tfrac{m(A\cap T^{-k}A)}{ m(A)^2},\ \ \text{and set }\ \
   u_A(\l):=\sum_{n=0}^\infty e^{-\l n} \tfrac{m(A\cap T^{-n}A)}{m(A)^2}. $$
The set $A\in\mathcal F$ is called a {\it moment} set for $T$ if for all $p\in\Bbb N$,
$$ \sum_{n=0}^\infty e^{-\l n}\int_AS_n(1_A)^pdm
\sim p!\,m(A)^{p+1} \frac{u_A(\l)^p}{\l}\ \ \text{ as }\l\to 0. $$

\begin{rem}
\sms {\bf a)} If $m(X)<\infty$ then (by the ergodic theorem)
every $A\in\B$ is a moment set.\sms
{\bf b)} Any conservative, ergodic, measure preserving transformation  with  moment
sets is rationally ergodic. Thus, for example, a  squashable conservative, ergodic,
measure preserving transformation (which is not rationally ergodic, see [A1])
has no moment sets.
\sms {\bf c)} We extend both
result and method from pointwise dual ergodic systems (as in [A1]) to
weakly pointwise dual ergodic situations. A similar approach was used in [T2]
to prove an arcsine-type limit theorem for pointwise dual ergodic maps.
We do not know if the latter result generalizes accordingly.
\end{rem}

\begin{theorem}[{\bf Moment set theorem}]\label{T_MomentSets}
\
Suppose that  $T$ is weakly  pointwise dual
ergodic, and that $A\in\mathcal F$ is a limited set.
 Then $A$ is a  moment set for $T$.
\end{theorem}

In view of Theorem 3.6.4 of [A1], the existence of limited sets established in
Proposition \ref{Prop_WPDE_local} above, and Karamata's Tauberian theorem (cf. p. 116 of [A1] or
Theorem 1.7.1 of [{BGT}]), Theorem \ref{T_DK} above is an immediate consequence of this result.
To prove the Moment set theorem, we need the following observation:

\begin{lem}[{\bf Convergence of Laplace transforms}]\label{Lem_CgeLaplace}\
Suppose that  $T$ is a c.e.m.p.t. on $(X,\mathcal B,m)$, weakly pointwise dual ergodic
with $A$ a limited set. Then
\begin{equation}\label{eq_a}
R_A(\l) := {1\over u_A(\l)}\sum_{n=0}^\infty e^{-\l n}\T^n 1_A
\underset{\l\searrow 0}{\overset{m_A}{\lra}} m(A),
\end{equation}
and
\begin{equation}\label{eq_b}
\varlimsup_{\l\searrow 0}\, R_A(\l) = m(A)\ \ \text{\rm a.e. on A},
\end{equation}
as well as
\begin{equation}\label{eq_cc}
\varlimsup_{\l\searrow 0}\, \sup_A\, R_A(\l)  = m(A)\ \ \text{\rm a.e. on A}.
\end{equation}
Also, each sequence $\l_i \searrow 0$ contains a subsequence
$\l'_i \searrow 0$ for which
\begin{equation}\label{eq_uhu1}
R_A(\l'_i) \underset{i \to \infty}{\lra} m(A)\ \ \text{\rm a.e. on A}.
\end{equation}
Moreover, every $B \in \mathcal F$, $B \subset A$, satisfies $a_n(A) \sim a_n(B)$
as $n \to \infty$, and hence, for $\l \searrow 0$,
\begin{equation}\label{eq_uhu2}
u_A(\l) \sim u_B(\l), \ \ \text{ and } \ \
\frac{R_A(\l)}{R_B(\l)} \to \frac{m(A)}{m(B)}\ \ \text{\rm a.e. on X}.
\end{equation}
In particular, for any sequence $(\l'_i)$ as in (\ref{eq_uhu1}), we also have
\begin{equation}\label{eq_uhu3}
R_B(\l'_i) \underset{i \to \infty}{\lra} m(B)\ \ \text{\rm a.e. on A}.
\end{equation}
\end{lem}

\pf \, Multiplying numerator and denominator by $(1-e^{-\l})$, we get
$$ 0 \le R_A(\l) = \frac{\sum_{n=0}^\infty e^{-\l n}\sum_{k=0}^{n-1}\T^k 1_A}{\sum_{n=0}^\infty e^{-\l n} a_n(A)}$$
for $\l>0$. As $A$ is a limited set, and $\T^k 1_A \le 1_X$, we therefore see that
(since $\sum_{k=0}^{\infty}\T^k 1_A = \infty$ a.e.),
\begin{equation}\label{eq_c}
\varlimsup_{\l\searrow 0}\, \sup_A\, R_A(\l)  \le m(A).
\end{equation}
On the other hand, monotone convergence ensures that
$\int R_A(\l)\, dm_A = m(A)$ for all $\l >0$. Together with
(\ref{eq_c}) this yields (\ref{eq_a}), and combining
(\ref{eq_c}) with (\ref{eq_a}) proves (\ref{eq_b}).
Together with (\ref{eq_c}) the latter gives (\ref{eq_cc}).
It is a standard fact from integration theory that
sequences which converge in probability contain a.e.
convergent subsequences, whence (\ref{eq_uhu1}).\\

Fix $B \in \mathcal F$, $B \subset A$.
We have $a_n(B)=m(B)^{-2} a_n(A) \int_B g_n(B)\, dm$
where $g_n(B) := a_n(A)^{-1} \sum_{k=0}^{n-1}\T^k 1_B {\overset{m_B}{\lra}} m(B)$
by weak pointwise dual ergodicity. Since $0 \leq T^k 1_B \leq T^k 1_A $
with $A$ a limited set, we see that
$\sup_n \sup_A g_n(B) < \infty$. Therefore,
$\int_B g_n(B)\, dm \to m(B)^2$, and hence $a_n(A) \sim a_n(B)$.
It is then immediate that $u_A(\l) \sim u_B(\l)$ since
$\sum_{n\geq0} a_n(A) = \infty$.
Now, expanding by $(1-e^{-\l})$ as above, we get
$$  \frac{R_A(\l)}{R_B(\l)} =
    \frac{u_B(\l)}{u_A(\l)}
    \frac{\sum_{n=0}^\infty e^{-\l n}\sum_{k=0}^{n-1}\T^k 1_A}
    {\sum_{n=0}^\infty e^{-\l n}\sum_{k=0}^{n-1}\T^k 1_B}
    \to 1\cdot\frac{m(A)}{m(B)}\ \ \text{\rm a.e. on X} $$
by Hurewicz' theorem (\ref{eq_Hu}) since
$\sum_{k=0}^{\infty}\T^k 1_A = \infty$ a.e. This proves
(\ref{eq_uhu2}), and (\ref{eq_uhu3}) follows at once.\; \qed\\

We are now ready for the\\

\pf {\it of the Moment set theorem}\;
{\bf (i)} We amend the argument given in the proof of Theorem\ 3.7.2 in [{A1}], using
the same combinatorial decomposition
\begin{equation}\label{eq_Diamond}
S_n(1_A)^p = \sum_{q=1}^p \g_p(q)\,a(q,n),
\end{equation}
where, for $n,p\in\Bbb N$, $a(p,n):X\to\Bbb Z$ is defined by
$a(0,n)(x):= 1$, and $a(p+1,n)(x):=\sum_{k=1}^n1_A(T^kx)a(p,n-k)(T^kx)$,
while $\g_1(q):=\d_{1,q}$, and $\g_{p+1}(q):=q(\g_p(q)+\g_p(q-1))$.
In particular, $\g_p(p)=p!$.\\

Proving that $A$ is a moment set reduces to showing that for $p\geq0$,
\begin{equation}\label{eq_Dontwash}
u_p(\l):=\sum_{n=0}^\infty e^{-\l n}\int_Aa(p,n)dm
\underset{\l\to 0}\sim m(A)^{p+1}{u_A(\l)^p\over\l}.
\end{equation}
This is because $p!\,u_p(\l)$, the $q=p$ term of the Laplace transform of the sum in (\ref{eq_Diamond}), dominates the $q<p$ terms. Indeed,
as we now check by induction on $p$,
\begin{equation}\label{eq_diamond}
u_p(\l)=O\({u_A(\l)^p\over\l}\)
\text{ as }\l\to 0\ \forall\ p \geq 0.
\end{equation}
For $p=0$ this is evident. More precisely, we have
\begin{equation}\label{eq_diamond2}
u_0(\l) = \frac{m(A)}{1-e^{-\l}} \sim \frac{m(A)}{\l} \
\text{ as }\l\to 0.
\end{equation}
To pass from $p-1$ to $p$,
use the recursive relation
\begin{align*}
u_p(\l)  &=\sum_{n=0}^\infty e^{-\l n}
\sum_{k=1}^n \int_{A} (1_A\, a(p-1,n-k))\circ T^kdm\\
&= \sum_{n=0}^\infty e^{-\l n}\sum_{k=1}^n\int_{A}\T^k1_{A}\,a(p-1,n-k)dm\\
&= \int_A\(\sum_{k=1}^\infty e^{-\l k}\T^k1_A\)
\(\sum_{n=0}^\infty e^{-\l n}a(p-1,n)\)dm,
\end{align*}
and combine it with Lemma \ref{Lem_CgeLaplace}. This proves (\ref{eq_diamond}) and
hence sufficiency of (\ref{eq_Dontwash}).\\

Since the $p=0$ case of (\ref{eq_Dontwash}) is trivially fulfilled,
we can establish (\ref{eq_Dontwash}) by proving
\begin{equation}\label{eq_1}
\liminf_{\l\to 0} {\l u_p(\l)\over u_A(\l)^p}\ge  m(A)^{p+1}\ \ \forall\
p\in\Bbb N
\end{equation}
\nobreak
and
\begin{equation}\label{eq_2}
\limsup_{\l\to 0} {\l u_p(\l)\over u_A(\l)^p}\le  m(A)^{p+1}\ \
\forall\ p\in\Bbb N.
\end{equation}\\

{\bf (ii)} Fix any $p\ge1$. To prove (\ref{eq_1}) by contradiction, suppose that it is violated.
Then there are some $\e>0$ and $\l_i \searrow 0$ such that
\begin{equation}\label{eq_ddagger}
\l_i\, u_p(\l_i) < (1-\e)^{p+2} \, m(A)^{p+1} u_A(\l_i)^{p}\ \;\;\;\;\; \text{for}\ i\ge1.
\end{equation}
Let $(\l'_i)$ be a subsequence of $(\l_i)$ as in Lemma \ref{Lem_CgeLaplace},
so that
\begin{equation}\label{eq_dddd}
R_B(\l'_i) \underset{i \to \infty}{\lra} m(B)\ \ \text{\rm a.e. on A}
\end{equation}
for all $B \in \mathcal F$, $B \subset A$.
We claim that there are measurable sets
$A=A_0 \supset A_1 \supset \ldots \supset A_p$ with $m(A_j)>(1-\e)\,m(A)$ and
$$  R_{A_j}(\l'_i) = \frac{1}{u_{A_j}(\l'_i)}\, \sum_{k \ge 1} e^{-\l'_i k}\, \T^k 1_{A_j} \ge (1-\e)\,m(A_j)\,\ \;\;\;\;\;
\text{on}\ A_{j+1}\ \text{for}\ i\ge \ell_j, $$
with $(\ell_j)_{j\ge1}$ increasing in ${\mathbb N}$.\\

To see this, start with $A_0=A$ and consider (\ref{eq_dddd})
with $B=A_0$. By Egorov's theorem, there is
some $A_1 \in \B \cap A_0$ with $m(A_1)>(1-\e)\,m(A)$ such that
this convergence is uniform on $A_1$. Therefore, there exists
a suitable $\ell_1$. If $A=A_0 \supset A_1 \supset \ldots \supset A_j$
have been constructed, consider (\ref{eq_dddd})
with $B=A_j$. By Egorov's theorem, there is
some $A_{j+1} \in \B \cap A_j$ with $m(A_{j+1})>(1-\e)\,m(A)$ such that
this convergence is uniform on $A_{j+1}$. Therefore, there exists
a suitable $\ell_{j+1}$.\\

We now find, for any $j\in\{0,1,\ldots,p-1\}$ and $i>\ell_j$, that
\begin{align*}
&\int_{A_j} \( \sum_{n=0}^\infty e^{-\l'_i n} a(p-j,n) \)\,dm \\
&=\sum_{n=0}^\infty e^{-\l'_i n}
\sum_{k=1}^n \int_{A_j} (1_A\, a(p-j-1,n-k))\circ T^kdm\\
&= \sum_{n=0}^\infty e^{-\l'_i n}\sum_{k=1}^n\int_{A}\T^k1_{A_j}\,a(p-j-1,n-k)dm\\
&= \int_{A}\(\sum_{k=1}^\infty e^{-\l'_i k}\T^k1_{A_j}\)
\(\sum_{n=0}^\infty e^{-\l'_i n}a(p-j-1,n)\)dm\\
&\ge \int_{A_{j+1}}\(\sum_{k=1}^\infty e^{-\l'_i k}\T^k1_{A_j}\)
\(\sum_{n=0}^\infty e^{-\l'_i n}a(p-j-1,n)\)dm\\
&> (1-\e)m(A)u_{A_j}(\l'_i)\int_{A_{j+1}}\(\sum_{n=0}^\infty e^{-\l'_i n}a(p-j-1,n)\)dm.
\end{align*}
Putting these together, we obtain that for $i>\ell_{p-1}$,
 \begin{align*}u_p(\l'_i) &>
((1-\e)m(A))^{p} \(\prod_{j=0}^{p-1} u_{A_j}(\l'_i)\)  \int_{A_{p}}\(\sum_{n=0}^\infty e^{-\l'_i n}a(0,n)\)dm\\ &
>((1-\e)m(A))^{p+1} \(\prod_{j=0}^{p-1} u_{A_j}(\l'_i)\) \frac{1}{\l'_i} ,\end{align*}
where the last step is immediate from $a(0,n)=1$.
Since $u_{A_{j}}(\l) \sim u_A(\l)$ by Lemma \ref{Lem_CgeLaplace},
this contradicts (\ref{eq_ddagger}), and thus establishes (\ref{eq_1}).\\

{\bf (iii)}  To prove (\ref{eq_2}), fix any $p\geq1$, and $\e>0$.
In view of (\ref{eq_cc}) there is some $\l'>0$ such that
$R_A(\l) < (1+\e) m(A)$ on A for $\l<\l'$. For such $\l$
we therefore find
\begin{align*}
u_p(\l) &=      \int_{A}\ u_A(\l) R_A(\l) ( \sum_{n=0}^\infty e^{-\l n} a(p-1,n) \)\,dm \\
        &\leq  (1+\e) m(A) u_A(\l) \cdot u_{p-1}(\l) \\
        &\vdots\\
        &\leq   ((1+\e) m(A) u_A(\l))^{p}  \cdot u_{0}(\l),
\end{align*}
and our claim is immediate from (\ref{eq_diamond2}).\ \qed\\

\section{Wandering rates, return sequences and tails of marginals}

\subsection*{\sl Wandering rates}

Suppose that  $(X,\B,m,T)$ is a c.e.m.p.t.
The {\it wandering rate} of the set $A\in\mathcal F$ is
the sequence given by $L_A(n):=m(\bigcup_{k=0}^nT^{-k}A)$, $n\ge1$.
Evidently,
$$A,B\in\mathcal F,\ A\subset B \Rightarrow L_A(n)\le L_B(n),$$
and
$$ \text{for}\ N\ge 1\ \text{fixed,}\ \; L_{\bigcup_{k=0}^NT^{-k}A}(n)=L_A(n+N)
\underset{n \to \infty}{\sim} L_A(n).$$
Wandering rates are expectations of truncated return times,
\begin{equation*}
L_A(n) = \int_A (\v_A \wedge n) dm.
\end{equation*}
Therefore, letting $c_A(\l):=\int_A(1-e^{-\l\v_A})dm$, $\l>0$,
for $A\in\mathcal F$, we have
$$c_A(\l) = (1-e^{-\l})\sum_{n=0}^\infty e^{-\l n} m(A\cap [\v_A>n])
\underset{\l\searrow 0}\sim \l^2\sum_{n=0}^\infty e^{-\l n} L_A(n).$$
Thus if $L_A(n)\sim L_B(n)$ as $n\to\infty$,
then $c_A(\l)\sim c_B(\l)$ as $\l\searrow 0$.
In fact, since $L_A(n+1)-L_A(n)\searrow 0$ for all $A\in\mathcal F$,
Korenblum's ratio\ Tauberian \ theorem ([{K}], see also Theorem 2.10.1 of [{BGT}])
shows that the converse is also true, so that
\begin{equation}\label{eq_IroningI}
\text{for}\ A,B\in \mathcal{F}: \ \;\;\;
L_A(n)\underset{n\to\infty}\sim L_B(n)\ \iff\ c_A(\l) \underset{\l\searrow 0}\sim c_B(\l).
\end{equation}

The set $A\in\mathcal F$ is said to have {\it minimal wandering rate } if
$L_B(n)\sim L_A(n)$ for all $B\in \mathcal{F}$, $B \subseteq A$.  In this case,
$\liminf_{n\to\infty}\tfrac{L_B(n)}{L_A(n)}\ge 1$ for all $B\in\mathcal F$.
Thus if $A, B\in\mathcal F$ both have minimal wandering rate,
then $L_B(n)\sim L_A(n)$, which defines the {\it wandering rate of} $T$,
$(L_T(n))_{n\ge 1}$ up to asymptotic equivalence. There are sufficient
conditions for $A\in\mathcal F$ to have minimal wandering rate.
By Proposition 3.2, Remark 3.6, and equation (2.3) of [{TZ}],
\begin{align}\label{eq_wr}
\text{if}\ ( \tfrac{\T_{A} (\v_A\wedge n)}{L_A(n)}  )_{n\ge 1}\
\text{is uniformly integrable,}\nonumber\\
\text{then}\ A\ \text{has minimal wandering rate.}
\end{align}

Also, uniform sets are known to have minimal wandering rate, provided that
the return sequence is regularly varying (Theorem 3.8.3 of [{A1}]).
In Theorem \ref{T_minWR} below we remove the latter condition.\\

Minimal wandering rates determine the return sequence $(a_n(T))_{n\ge1}$
of a weakly pointwise dual ergodic system $(X,\B,m,T)$ by means of the
{\it asymptotic renewal equation}. Assuming w.l.o.g. that $a_n(T)=\sum_{j=0}^{n-1} u_n(T)$
with $u_n(T)\ge0$, we let
$$u_T(\l) := \sum_{n=0}^\infty e^{-\l n} u_n(T),\ \;\;\;\; \text{for}\ \l>0.$$
As a consequence of (\ref{eq_RS_A}) we have $u_T(\l) \sim u_A(\l)$ as $\l \searrow 0$
for all $A\in \mathcal L(T)$ (with $u_A(\l)$ as in DK-section).

\begin{theorem}[{\bf Minimal wandering rates and the asymptotic renewal equation}]\label{T_minWR}
Suppose that  $T$ is weakly  pointwise dual ergodic.\medskip

\noindent {\bf(i)}\ If $A\in\mathcal{F}$ has
minimal wandering rate, then it satisfies the asymptotic renewal equation
$$c_A(\l) \underset{\l\to 0+}\sim \frac1{u_T(\l)}.$$\medskip
\noindent {\bf(ii)}\ Uniform sets have minimal wandering rates.
\end{theorem}

\pf {\bf (i)}\ By Proposition \ref{Prop_WPDE_local} there is some limited set $B\in \mathcal{F}_A$.
In view of statement a) in Lemma \ref{Lem_CgeLaplace}, any sequence decreasing to $0$ contains a subsequence
$(\l_j)_{j\ge 1}$ along which
$$\frac{1}{u_B(\l_j)} \sum_{n=0}^\infty e^{-\l_j n}\T^n 1_B
\underset{j \to \infty}\lra m(B)\ \;\;\;\; \text{a.e. on}\ B.$$
Egorov's theorem then gives us some $C\in \mathcal{F}_B$ on which this
convergence is in fact uniform, so that
$$\int_C (1-e^{-\l_j\v_C})\sum_{n=0}^\infty e^{-\l_j n}\T^n 1_B dm
\underset{j \to \infty}\sim m(B)\, u_B(\l_j)\, c_C(\l_j).$$
On the other hand,  Lemma 3.8.4 in [A1] shows that
$$\int_C (1-e^{-\l\v_C})\sum_{n=0}^\infty e^{-\l n}\T^n 1_B dm
=\sum_{n=0}^\infty e^{-\l n}\int_{C_n} 1_B dm \underset{\l \searrow 0}\lra 1,$$
where $C_0:=C$ and $C_n:=T^{-n}C\setminus\bigcup_{k=0}^{n-1}T^{-k}C$ for $n\ge 1$.
Together, these prove $c_C(\l_j) \sim 1/{u_B(\l_j)}$ as $j \to \infty$.
But as $A$ has minimal wandering rate, (\ref{eq_IroningI}) ensures that
$c_A(\l)\sim c_C(\l)$, and we end up with
$c_A(\l_j)\sim 1/{u_B(\l_j)}\sim 1/{u_T(\l_j)}$. Our claim follows
since this can be done inside any sequence of $\l$'s decreasing to $0$.\\

{\bf (ii)}\ Suppose that $A\in\mathcal U(T)$ is uniform for $f\in L_+^1(m)$,
with return sequence $(a_n)_{n\ge 1}$.
Then, this is also true for all $B\in \mathcal{F}_A$. Thus, by the asymptotic
renewal equation for uniform sets (cf. 3.8.6 of [A1]),
$$c_B(\l) \underset{\l \searrow 0}\sim \tfrac1{u_T(\l)}\ \;\;\;\; \text{for}\ B\in \mathcal{F}_A,$$
where $u_T(\l)$ does not depend on $B$.
In particular, $c_B(\l)\sim c_A(\l)$ for all $B\in \mathcal{F}_A$,
which, due to (\ref{eq_IroningI}), shows that $A$ has minimal wandering rate.\ \qed\\

\subsection*{Proof of Theorem \ref{T_NormalizationID}}

Let $(X,\B,m,T)$ be the Kakutani tower of $(\Om,\mathcal A ,P,S,\v)$, then
$$L_\Om(n) = \Bbb E(\v\wedge n) \sim \tfrac{n}{(1+\g)A(n)},$$
whence by Theorem 3.8.1 of [{A1}], for large $n$,
$$a(n)=\sum_{k=1}^nm(\Om\cap Y^{-k}\Om)\ge\tfrac{n}{2L_\Om(n)}>\tfrac12 A(n).$$
Thus, for all $\e>0$,
$$\tfrac{n\vartheta_\mu(\e a(n))}{a(n)}<\tfrac{2n\vartheta_\mu((\e/2) A(n))}{A(n)}\underset{n\to\infty}\lra 0.$$
Now, by Theorem \ref{Thm_WPDE_Thetamu} (a), $T$ is weakly $\pde$ with return sequence $a(n)$, and
in view of (\ref{eq_Wheelchair}) and (\ref{eq_wr}), $\Om$ has minimal wandering rate.
According to the asymptotic renewal equation of Theorem \ref{T_minWR},
$c_\Om(\l)\underset{\l\to 0+}\sim\tfrac1{u(\l)}$ whence by Karamata's theorem
$a(n)\sim A(n)$.\ \ \ \qed

\section{The one-sided law of the iterated logarithm}

The $\g=1$ version of the law of the iterated logarithm
follows immediately from the previous results.

\subsection*{Proof of Theorem \ref{T_SLT}(b)}
It has already been pointed out in [AD1] that (\ref{eq_LIL1})
holds for positive stationary processes satisfying a
weak law of large numbers provided that the corresponding
infinite measure preserving Kakutani tower is weakly
pointwise dual ergodic. The latter is immediate from
Theorem \ref{Thm_WPDE_Thetamu}.
\qed\\

We now prove Theorem \ref{T_osLIL} by applying [{AD2}].

\subsection*{Proof of Theorem \ref{T_osLIL}}
We first show that under the present assumptions
\begin{equation}\label{eq_Stopsign}
\sum_{n=1}^\infty\tfrac{\phi_-(n)}n<\infty.
\end{equation}
To see this, note that $\phi_-(a(a(n)))\le \tfrac{a(n)}{n}$ for large $n$. Let $b$ be
asymptotically inverse to $a$ in that $b(a(n))\sim a(b(n))\sim n$,
then $b$ is $\tfrac1\g$-regularly varying, and for large $N:=a(a(n))$ we have
$$\phi_-(N)=\phi_-(a(a(n)))\le \tfrac{a(n)}{n}=\tfrac{b(N)}{b(b(N))}=\tfrac1{c(N)}$$
 where $c(N):=\tfrac{b(b(N))}{b(N)}$ is $(\tfrac1{\g^2}-\tfrac1\g)$-regularly varying.
 Since $Nc(N)$ is $(\tfrac1{\g^2}-\tfrac1\g+1)$-regularly varying
we indeed get $\sum_{N=1}^\infty\tfrac1{Nc(N)}<\infty$.

As an immediate consequence of (\ref{eq_Stopsign}),
$(\Om,\mathcal A ,P,S,\a)$ is {\it strongly mixing from below} as defined in [{AD2}].\\

Let $(X,\B,m,T)$ be the Kakutani tower of $(\Om,\mathcal A ,P,S,\v)$.
Part c) of Theorem \ref{Thm_WPDE_Thetamu}, with $\mu = P$,
ensures that $T$ is pointwise dual ergodic and $\Om\in\mathcal{DK}(T)$.\\

The assumptions of  Theorem 4 in [{AD2}] are now satisfied. Hence,
\begin{align*}
&\sum_{n=1}^\infty \tfrac 1ne^{-\beta\tau(n)
}<\infty\ \forall\ \beta>1\ \implies\\ &  \varlimsup_{n\to\infty}\tfrac 1{a(n/\tau(n))\tau(n)}
S_n(f)
\le K_\g  \int_Xf\,d\mu\ \text{ a.e.}\ \forall\ f\in L_+^1
\end{align*}
and
\begin{align*}
&\sum_{n=1}^\infty\tfrac 1ne^{-r\tau(n)
}=\infty\ \forall\ r<1\ \implies\\ &  \varlimsup_{n\to\infty}\tfrac 1{a(n/\tau(n))\tau(n)}
S_n(f)
\ge K_\g  \int_Xf\,d\mu\ \text{ a.e.}\ \forall\ f\in L_+^1.
\end{align*}

Using the inversion technique in \S5 of [{AD2}], statements
(a) and (b) of Theorem \ref{T_osLIL} follow, and (c)
is a consequence of (a) and (b).\qed\\

\section{Interval  maps}
 A {\it piecewise monotonic} ({\it increasing}) {\it map of the interval} is a triple
$(\Om,S,\a)$ where $\Om$ is a bounded interval, $\a$ is a finite or countable
generating partition (mod $m:=$ Lebesgue measure) of $\Om$ into open intervals,
and $S:\Om\to \Om$ is a map such that the restriction $S:A\to SA$ is an {(increasing)} homeomorphism
for each $A\in\a$ so that both $S:A\to SA$  and its inverse $v_A:SA\to A$ are absolutely continuous.

In this case, each iterate $(\Om,S^k,\a_k)$, $k \ge 1$, is also piecewise monotonic (increasing),
where $\a_k :=\bigvee_{i=0}^{k-1} S^{-i} \a$. Generalizing the above, we let,
for $A \in \a_k$, $v_A$ denote the inverse of $S^{k}:A\to S^{k}A$, so that the
transfer operator of $S$ (with respect to $m$) satisfies
$$
\widehat S ^k f=\sum_{A\in\a_k}1_{S^kA}v_A'\,\cdot\,(f\circ v_A),
\textrm{ where }
v_A':=\frac{dm\circ v_A}{dm}.
$$

Consider
the following properties for a piecewise monotonic map of the interval $(\Om,S,\a)$:
\begin{enumerate}
\item[(A)] {\it Adler's condition\/}: for all $A\in\a$, $S|_A$ extends to a
$C^2$ diffeomorphism on $\overline{A}$, and
$S''/(S')^2$ is bounded on $\Om$.
\item[(B)] {\it Big images\/}: $\min_{A\in\a}m(SA)>0$.
\item[(R)] {\it Rychlik's condition\/}: $\sum_{A\in\a}\|1_{SA}v_A'\|_{\widehat{BV}}=:\mathcal R<\infty.$
\item[(U)] {\it Uniform expansion\/}: $\inf|S'|>1$.
\end{enumerate}
Recall that (A) ensures $\tfrac{v_A^{\pprime}}{v_A^{\prime}}\le M < \infty$,  whence
$v_A^{\prime}=e^{\pm M}\tfrac{m(A)}{m(SA)}$ on $SA$ for all  $A\in\a$.
In (R), the space $\widehat{BV}$ is the subspace of those functions in $L^\infty(m)$ with a version in $BV$, the space of functions of bounded variation. The norm $\|\cdot\|_{\widehat{BV}}$ is defined by
$$
\|f\|_{\widehat{BV}}:=\|f\|_\infty+\widehat{\bigvee_\Om} f,\;\textrm{ where }\widehat{\bigvee_\Om}
f:=\inf\{\bigvee_\Om(f^\ast):f^\ast=f\ m\text{-a.e.}\}.
$$
 Piecewise monotonic maps $(\Om,S,\a)$ of the interval with
properties (P$_1$),...,(P$_N$) will be called
{\em P$_1$...P$_N$ maps} (eg  ABU, RU maps).\\

\begin{lem}\label{Prop_ABU_RU}
Any ABU map is an RU map.
\end{lem}

\pf.  This is similar to Proposition 2 of [{Z1}].
For any piecewise monotonic map $(\Om,S,\a)$, (A)\&(B) imply (R).
Indeed,
\begin{align*}
\sum_{A\in\a}\|1_{SA}v_A'\|_{\widehat{BV}}& \le 3\|v_A'\|_\infty+\bigvee_{SA}(v_A')\\ &\overset{\textrm{(A)}}\le
\sum_{A\in\a}(3e^M\tfrac{m(A)}{m(SA)}+\int_{SA}|v_A^{\pprime}|dm)\\ &\le
\sum_{A\in\a}(3e^M\tfrac{m(A)}{m(SA)}+M\int_{SA}v_A^{\prime}dm)\\ &\le
\sum_{A\in\a}(3e^M\tfrac{m(A)}{m(SA)}+Me^Mm(A))\\ &\overset{\textrm{(B)}}\le M'\sum_{A\in\a}m(A)=M'.\qed
\end{align*}

\subsection*{\sl Ergodic properties of Rychlik's maps}
Suppose that $(\Om,S,\a)$ is a RU map, then, according to [{R}],
\bul $(\Om,\B,m,S,\a)$ is a fibred system where $m$ is Lebesgue measure on $\Om$ and $\B$ denotes the Borel $\sigma$-field,
\bul the ergodic decomposition of $(\Om,\B,m,S)$ is finite,
\bul to each ergodic component there corresponds an absolutely continuous invariant probability, with density in $BV$ and with respect to  which $S$ is isomorphic to
the product of a finite permutation and a mixing RU map.

If  $S$ is {\it weakly mixing} (with respect to $m$ in the sense that
$f:\Om\to\Bbb S^1$ measurable, $f\circ S=\l f$ a.e.
where $\l\in\Bbb S^1$ implies $f$ constant), then there are
constants $K>0 ,\, \th\in (0,1)$ such that
$$
\left\| \widehat S ^n f - \left(\int_\Om f dm \right) h \right\|_{\widehat{BV}} \leq K \th^n \|f\|_{\widehat{BV}}
$$
\f where $h$ is the unique $T$-invariant probability density.
In this case, let $dP:=hdm$ and $\mu:=m|_{[h>0]}$.
Then [{AN}] shows that the probability preserving fibred system
\begin{equation}\label{eq_ExpThetaMuMix}
(\Om,\B,P,S,\a) \text{ is exponentially } \vartheta_\mu \text{-mixing.}
\end{equation}

We next observe that $\psi^*(N)<\infty$ already implies cf-mixing
in the present context, provided $h$ has a positive lower bound.
(This shows that for such ABU maps, the conclusions of Theorem \ref{T_osLIL}
already follow from earlier results for cf-mixing systems.)

\begin{propn}[{\bf cf-mixing ABU maps}]\label{Prop_CFmixing}
Let $(\Om,S,\a)$ be a weakly mixing ABU map with invariant density $h$ bounded away from $0$.
If $\exists\ N\ge 1\ \st\ \ \psi^*(N)<\infty$, then $(\Om,S,\a)$ is cf-mixing.
\end{propn}

\pf.\;  Suppose that $\eta\in (0,1)$ satisfies $h=\eta^{\pm 1}$, which we use to
abbreviate $\eta \le h \le \eta^{-1}$. A standard argument then shows that
$\sup_{n\ge 1}\sup_\Om|(S^{n})^{\pprime}|/((S^{n})^{\prime})^2<\infty$, and we can also assume that
$$v_A'=\eta^{\pm 1}\tfrac{m(A)}{m(S^kA)}\ \; \text{ on } S^kA\ \;\;\;\text{ for all }  A\in\a_A, k\ge 1.$$
Let $\widehat S_P$ be the transfer operator with respect to the absolutely continuous invariant probability $P$,
then $\widehat S_P f = \widehat S(hf)/h$, and therefore
$$\widehat S_P^nf = \eta^{\pm 2}\widehat S^nf\ \;\;\;\text{ for all }   n\ge 1,\ f\in L_+^\infty.$$
We now show that
\begin{equation}\label{eq_ImgEsti}
m(S^kA)\ge \D\ \;\;\; \text{ for all } A\in\a_k,\ k\ge 1,
\;\;\;\text{ where } \D:=\tfrac{\eta^6}{\psi^*(N)}.
\end{equation}
To this end, let $B\subset \Om$ be measurable with $m(B)>0$. Then
\begin{align*}
\frac{P(A)}{m(S^kA)} &\le\eta^{-1}\frac{m(A)}{m(S^kA)}\\ &
\le \eta^{-2}\frac1{m(S^{-N}B)}\int_{S^{-N}B}1_{S^kA}v_A' dm\\ &=
 \eta^{-2}\frac1{m(S^{-N}B)}\int_{S^{-N}B} \widehat S^k1_A dm\\ &\le
\eta^{-6}\frac1{P(S^{-N}B)}\int_{S^{-N}B} \widehat S_P^k1_A dP\\ &=
\eta^{-6}\frac1{P(B)}P(A\cap S^{-N+k}B)\\ &\le
\eta^{-6}\psi^*(N)P(A)
\end{align*}
whence (\ref{eq_ImgEsti}), as claimed.

To complete the proof of the proposition, we can then
proceed as in the proof of Theorem 1(b) in [AN].\ \ \qed

\

By exponential $\vartheta_\mu$-mixing (\ref{eq_ExpThetaMuMix}),
Theorem \ref{T_SLT} implies the general

\begin{propn}[{\bf Stable limit theorem for RU maps}]\label{T_SLT_ABU}
Suppose that $(\Om,S,\a)$ is a  weakly mixing RU map with absolutely continuous
invariant probability $dP=hdm$. Let $\v:\Om\to\Bbb N$ be $\a$-measurable,
and denote $\v_n:= \sum_{k=0}^{n-1} \v \circ S^k$.

If $a(n):=\sum_{k=1}^nP([\v_k\le n])$ is $\g$-regularly varying for $\g\in (0,1]$,
with inverse $b$, then
$$\frac{\v_n}{b(n)}\overset{\frak d}\lra\ Z_\g.$$
\end{propn}

\begin{rem}
For the subfamily of those RU-maps $S$ which satisfy (A) plus the
\emph{finite image condition (F)} which requires $\{SA:A\in\a\}$
to be finite, more general stable limit theorems (for observables
$\v$ which need not have constant sign) follow from [ADSZ] (see the
end of Section 5 there). These \emph{AFU maps} occur as induced maps of the
infinite measure preserving \emph{AFN maps} studied in [Z1],[Z2]
(generalizing [T1]). The final subsections below illustrate that the
present results allow us to analyse, via
weak pointwise dual ergodicity, infinite measure
preserving interval maps more general than those of [Z2].
\end{rem}

\

\subsection*{\sl The asymptotic type}
Next, we turn to the asymptotic identification, via
Theorem \ref{T_NormalizationID}, of the normalizing constants $a(n)$ in this setup.

\begin{propn}[{\bf Asymptotic type of ABU maps}]\label{Prop_amschluss}
Suppose that $(\Om,S,\a)$ is a weakly mixing ABU map
with absolutely continuous invariant probability $dP=hdm$.

Suppose that $\v:\Om\to\Bbb N$ is $\a$-measurable and satisfies
\begin{equation}\label{eq_aaaa}
\frac{P([\v\ge n])}{m([\v\ge n])}\underset{n\to\infty}\lra c\in (0,\infty),
\end{equation}
as well as
\begin{equation}\label{eq_bbbb}
\int_\Om \v\wedge n \, dm \underset{n\to\infty}\sim \frac{n}{\Gamma(2-\gamma)\Gamma(1+\gamma) A(n)},
\end{equation}
where
$A(t)$ is strictly increasing and regularly varying with index $\gamma\in (0,1]$. Then
$$a(n):=\sum_{k=1}^n P([ \v_k \le n]) \ \ \underset{n\to\infty}\sim\ \  c^{-1}A(n).$$
\end{propn}

The main point is condition (\ref{eq_Wheelchair}) of Theorem \ref{T_NormalizationID}.

\begin{lem}\label{L_amschluss}
Suppose that $(\Om,S,\a)$ is a weakly mixing ABU map with
absolutely continuous invariant probability $dP=hdm$.
Suppose that $\v:\Om\to\Bbb N$ is $\a$-measurable and satisfies
\begin{equation*} 
\int_\Om \v\wedge n \, dm = O \left( \int_\Om \v\wedge n \, dP  \right)
\ \ \ \text{ as } n\to \infty.
\end{equation*}
Then there is some $\Psi \in L^1(P)$ such that
$$\widehat S_P(\v\wedge n) \le \Psi \int_\Om \v\wedge n \, dP \ \ \text{ a.e. }
\ \ \ \text{ for all } n\ge 1.$$
\end{lem}

\pf. We first record a corresponding
statement w.r.t. Lebesgue measure $m$,
\begin{equation}\label{eq_iiiprime}
\exists\ \tilde{M}>0\ \text{  so that }\ \
\widehat S(\v\wedge n) \le \tilde{M}\int_\Om (\v\wedge n) dm\ \forall\ n\ge 1.
\end{equation}
Letting $F_n:=\sum_{A\in\a}(\v(A)\wedge n){m(A)}1_{SA}$
and $M':=(\inf_{A\in\a} m(SA))^{-1}$, we have indeed
\begin{align*}
\widehat S(\v\wedge n) &= \sum_{A\in\a}(\v(A)\wedge n)\widehat S 1_A
  = \sum_{A\in\a}(\v(A)\wedge n) v'_A1_{SA}\\
&\le e^{M}\sum_{A\in\a}(\v(A)\wedge n)\tfrac{m(A)}{m(SA)}1_{SA}
  \le M'e^{M}  F_n.
\end{align*}
But
$\|F_n\|_\infty \le \sum_{A\in\a}(\v(A)\wedge n){m(A)} = \int (\v\wedge n) dm$,
whence (\ref{eq_iiiprime}). To deduce (\ref{eq_Wheelchair}), note that
\begin{align*}
\widehat S_P(\v\wedge n) &=1_{[h>0]}\tfrac1h\,\widehat S (h (\v\wedge n))\\
&\le 1_{[h>0]}\|h\|_\infty \tfrac1h\,\widehat S(\v\wedge n)\\
&\le 1_{[h>0]}\|h\|_\infty \tfrac1h\, \tilde{M}\int_\Om (\v\wedge n) dm\ \ \ \text{ by (\ref{eq_iiiprime})}\\
&\sim 1_{[h>0]}\|h\|_\infty \tfrac1h\, c \tilde{M}\int_\Om (\v\wedge n) dP\\
&=\Psi \int_\Om (\v\wedge n) dP,
\end{align*}
where $\Psi:=1_{[h>0]}\|h\|_\infty \tfrac1h cM\in L^1(P)$ since $h\in L^1(m)$.\ \ \ \qed

\

\pf of Proposition \ref{Prop_amschluss}.  We are going to verify the conditions of Theorem \ref{T_NormalizationID}.
Note first that adaptedness follows from the other two
by exponential $\vartheta_\mu$-mixing (\ref{eq_ExpThetaMuMix}). Next,
condition (\ref{eq_NormIDcondition3}) is immediate from (\ref{eq_aaaa}), as
\begin{align*}
\int_\Om\v\wedge n\,dP & =\sum_{k=1}^nP([\v\ge k])
\overset{\text{ (\ref{eq_aaaa})}}{\underset{k\to\infty}\sim}\ c\sum_{k=1}^nm([\v\ge k])
=c\int_\Om\v\wedge n\,dm.
\end{align*}
To check condition (\ref{eq_Wheelchair}), use this and the previous lemma.\ \ \ \qed

 \

\subsection*{\sl The common image property }
Typically, for interval maps, one will first
obtain information on $[\v\ge n]$ in terms of Lebesgue measure $m$. This needs to be
combined with an analysis of $h$ to yield information on $P([\v\ge n])$, and hence
on $\int \v\wedge n \, dP = \sum_{k=1}^n P([\v\ge k])$.
Here we discuss simple sufficient conditions which allow us to validate property
(\ref{eq_aaaa}) of Proposition \ref{Prop_amschluss} in this way.

Consider a piecewise increasing map $(\Om,S,\a)$, with
$\Om=[\om_l,\om_r]$.
We shall say that $(\Om,S,\a)$ has the {\it common image property} if
$\bigcap_{A\in\a}SA=(\om_l,\om_l+z_S)$ where $z_S>0$. Evidently, this entails the
big image property (B). Moreover, we find:

\begin{lem}\label{Lem_LowBdH}
Suppose that
$(\Om,S,\a)$ is  a  piecewise increasing AU map with the
 common image property and an absolutely continuous invariant probability $dP=h\,dm$
 on $\Om=[\om_l,\om_r]$, then
 \begin{equation}\label{eq_999}
 \underset{[\om_l,\om_l+z_S]}{\text{\rm essinf}}\,h>0.
 \end{equation}
 Moreover, $S$ is weakly mixing.
\end{lem}

\pf. \; We assume w.l.o.g. that $\Om=[0,1]$.
Fix a version $h\in BV$ of the invariant density and set
$$ \mathcal J:=\{J\subset {[0,1]}: J \text{ is a nonempty open interval with } \underset{J}\inf\,h>0\}.$$
It is clear that $\mathcal J \neq \emptyset$. We need to show that
$(0,z_S)\in\mathcal J$. Observe first that
\begin{equation}\label{eq_99}
\text{there exist } J\in {\mathcal J}  \text{ and } A\in\a \text{ so that } J\cap\bdy A\ne\emptyset.
\end{equation}
To see this, suppose otherwise i.e. that
$\forall\ J\in\mathcal J,\ \exists\ A_J\in\a :\ J\subset A_J$.
Then $J\in\mathcal J$ implies $SJ\in\mathcal J$
since for $x\in SJ\subset SA_J$,
$$ h(x)\ge v_{A_J}'(x)h(v_{A_J}x)\ge \text{\rm const}\cdot m({A_J})\,\underset{J}{\inf}\,h>0.$$
But then, for each $k\ge 1$, $S^kJ\subset A_k$ for some $A_k\in\a$, an impossibility since this
entails $m(S^kJ)\ge\l^k m(J)\to\infty$.

Due to (\ref{eq_99}), there are $J\in\mathcal J$ and $A=(u,v)\in\a$ such that $u\in J$.
Set $J_0:=A \cap J=(u,w)$ with $u<w$. It follows as above that
$I_0:=SJ_0\in\mathcal J$, and the common image property implies
$I_0=(0,c)$ for some $c>0$.

Note then that there exist some $J'\in\mathcal J$ and $A'\in\a$ such that
$J'\supset A'$: Unless $I_0$ contains some $A'$, it is contained in a specific $A'\in\a$,
and by the special structure of our map there is some
$k\ge 1$ for which $S^kI\subset A' \subset S^{k+1}I$.
By the argument proving (\ref{eq_99}) we have $J':=S^{k+1}I \in \mathcal J$.

But then $(0,z_S)\subset SA'\in\mathcal J$ as required.

Finally, in view of Lemma \ref{Prop_ABU_RU} and [R], $S$ has only finitely
many ergodic acims, and these have densities $h_i \in BV$, which can be chosen to be
lower semicontinuous, so that the sets $[h_i>0]$ are open and pairwise disjoint.
However, by the above, each $[h_i>0]$ contains $(0,z_S)$. Hence $h$ and $P$ are unique, meaning
that $S$ is ergodic. Moreover, the structural results of [R] also show that there
is a finite tail decomposition $h=\sum_{j=0}^{p-1} g_j$ with
$g_j \in BV$ and the $[g_j>0]$ open and pairwise disjoint, such that
$S [g_j>0] = [g_l>0]$ a.e., $l=j+1 \mod p$, and $S$ is weakly mixing iff $p=1$.
Bounded variation of the $g_j$ together with (\ref{eq_999}) now implies that
(after renumbering the $g_j$ if necessary) there is some $y>0$ for which
$(0,y) \subseteq [g_0>0]$. However, there is at least one cylinder $A=(a,b) \in \alpha$
with $a < y$, and then $(a,c) := [g_0>0] \cap A$ is nonempty and open.
Due to the common image property, $S(a,c) \subseteq [g_1>0]$ has nonempty
open intersection with $(0,y)$. Hence $[g_0>0] = [g_1>0]$, as these sets overlap.
Therefore $p=1$.
\qed\\

This immediately allows us to deal with situations in which
$\v$ only diverges at $0$.

\begin{exple}\label{expleintervalmaps1}
Suppose that $([0,1],S,\a)$ is a pcw increasing AU map with the
common image property and absolutely continuous invariant probability $dP=hdm$.
Suppose that $\v:[0,1]\to\Bbb N$ is $\a$-measurable and satisfies
$\o{[\v\ge n]} = [0,y_n]$, where
$$ y_n \sim \frac{1}{\Gamma(1-\gamma) \Gamma(1+\gamma) A(n)}  $$
with $A$ strictly increasing and regularly varying of index $\gamma\in (0,1]$.
Then $a(n) \sim A(n)$ as $n\to\infty$.

Indeed, we need only check condition (\ref{eq_aaaa}) of Proposition \ref{Prop_amschluss}.
Fixing a version $h\in BV$ of the invariant density, Lemma \ref{Lem_LowBdH} shows that
$\lim_{x\to 0+}h(x)=:h(0^+)>0$. Whence
$P([\v\ge n])=\int_{[\v\ge n]}h\,dm\ \sim\ h(0^+)\,m([\v\ge n])$,
and our claim follows since, by Karamata's theorem,
$$ \int_\Om (\v\wedge n)dm \underset{n\to\infty}\sim\frac{n}{\Gamma(2-\gamma)\Gamma(1+\gamma) A(n)}.$$
\end{exple}

Next, we record a little preparatory observation which will enable us to also
study functions $\v$ which diverge at countably many points.
(This will be the case for the return time functions of the null-recurrent
maps studied in the final subsection below.)

\begin{lem}\label{Lem_ManyBits}
Let $h:\Om \to [0,\infty)$ have right-hand limits $h(x^+)$ everywhere. Let
$x_j,y_{j,n} \in \Om$, $j,n \geq 0$, be such that for each $n$ the sets
$(x_j,x_j+y_{j,n})$ are pairwise disjoint, and suppose that there are
$s_j \in [0,\infty)$ with $\sum_{j \ge 0} s_j < \infty$, and $q_n \searrow 0$
for which $ {y_{j,n}}/{q_n} \to s_j$  as $n \to \infty$, uniformly in $j$.
If $s_j h(x_j^+)>0$ for some $j$, then
 \begin{equation}\label{eq_99999}
 \sum_{j \ge 0} \int_{(x_j,x_j+y_{j,n})} h(x)\,dx \sim
 \left(\sum_{j \ge 0} s_j h(x_j^+) \right) \, q_n \ \ \text{ as } n \to \infty.
 \end{equation}
\end{lem}

 \pf. Assume w.l.o.g. that $s_0 h(x_0^+)>0$, and take any $\epsilon>0$.
 Choose $n_1$ so large that $y_{j,n} \le e^\epsilon s_j q_n$
 for $n \ge n_1$ and all $j$. Take $J \ge 1$ so large that
 $(\sup h) \sum_{j >J} s_j  < \epsilon \sum_{j \ge 0} s_j h(x_j^+)$.
 Next, there is some $n_2 \ge n_1$ such that for all $n \ge n_2$ and all $j \le J$,
 $$ \sup_{(x_j,x_j+e^\epsilon s_j q_n)} h \le
    e^\epsilon h(x_j^+) + \frac{\epsilon s_0 h(x_0^+)}{J+1}$$
 ($h(x_j^+)$ need not be positive, but $h(x_0^+)$ is).
 Then, for all $n \ge n_2$,
 \begin{align*}
  \sum_{j \ge 0} \int_{(x_j,x_j+y_{j,n})} h(x)\,dx
  &\le
  e^\epsilon q_n \sum_{j \ge 0} s_j \sup_{(x_j,x_j+e^\epsilon s_j q_n)} h \\
  &\le
  e^\epsilon q_n  \left( e^\epsilon \sum_{j = 0}^J s_j h(x_j^+)
  + \epsilon s_0 h(x_0^+) + \sup h \sum_{j > J} s_j  \right)   \\
  &\le
  e^\epsilon (e^\epsilon + 2 \epsilon) \left(\sum_{j \ge 0} s_j h(x_j^+) \right) q_n.
 \end{align*}
 As $\epsilon >0$ was arbitrary, this proves one half of our claim.
 The corresponding estimate from below follows by similar
 but even simpler arguments which we omit.
 \qed\\

We can now go beyond the scenario of
Example \ref{expleintervalmaps1} above.
Situations of the following type naturally occur in the study of interval maps
with neutral fixed points, see below.

\begin{propn}[{\bf $\v$ with countably many singularities}]\label{Prop_ganzamschluss}
Let $(\Om,S,\a)$ be a weakly mixing  pcw increasing ABU map with
absolutely continuous invariant probability $dP=h\,dm$.
Suppose that $\v:\Om\to\Bbb N$ is $\a$-measurable and such that for $n$ sufficiently large,
$[\v\ge n]$ is a countable disjoint union of intervals $(x_j,x_j+y_{j,n})$
satisfying the assumptions of Lemma \ref{Lem_ManyBits}, where
$$ q_n \underset{n\to\infty}\sim \frac{1}{\Gamma(1-\gamma) \Gamma(1+\gamma) A(n)}$$
with $A$ strictly increasing and regularly varying of index $\gamma\in (0,1]$.
Then
$$a(n):=\sum_{k=1}^n P([\v_k\le n])\ \ \underset{n\to\infty}
\sim\ \ \left(\sum_{j \ge 0} s_j h(x_j^+) \right) \, A(n).$$
\end{propn}

\pf.  It is clear that
$ m([\v\ge n]) \sim (\sum_{j \ge 0} s_j ) \, q_n$.
According to Lemma \ref{Lem_ManyBits},
$$ P([\v\ge n]) = \sum_{j \ge 0} \int_{(x_j,x_j+y_{j,n})} h(x)\,dx \sim
   \left(\sum_{j \ge 0} s_j h(x_j^+) \right) \, q_n, $$
so that by Karamata's theorem
$$ \int_\Om (\v\wedge n)dP \underset{n\to\infty} \sim
    \frac{\left(\sum_{j \ge 0} s_j h(x_j^+) \right) n}{\Gamma(2-\gamma)\Gamma(1+\gamma) A(n)},$$
and Proposition \ref{Prop_amschluss} applies.\ \ \qed

 \

\subsection*{\sl Some infinite measure preserving interval maps }
We conclude with a class of infinite measure preserving
interval maps $T$ with indifferent fixed point,
which induce probability preserving maps $S$ of the
above type. These $T$ do not, in general, belong to the family of
AFN-maps studied in [Z2].

\begin{propn}[{\bf Maps with indifferent fixed points}]\label{Prop_ganzganzamschluss}
Let $(X,T,\b)$ be a pcw increasing A map on $X=[\eta_l,\eta_r]$ with the common image property which
satisfies $\inf_{(\eta_l+\epsilon,\eta_r)} T'>1$ for every $\epsilon \in (0,\eta_r-\eta_l)$.
Assume that $T$ possesses  a leftmost cylinder $B_*=(\eta_l,\xi)$, and that
$z_T := \inf_{B \in \b} m(TB)$ satisfies $z_T>\xi-\eta_l$.
Suppose that $T$ is convex near $\eta_l$, and satisfies, for some $\gamma \in (0,1)$,
\begin{equation}\label{eq_indiff}
T(\eta_l+x) \sim \eta_l+x + \kappa x^{1+1/\gamma} + o(x^{1+1/\gamma})\, \ \ \text{ as }\ x \searrow 0.
\end{equation}

Then $T$ is conservative ergodic with an infinite acim $m_T = h_T \, dm$,
with $h_T$ bounded on each $(\eta_l+\epsilon,\eta_r)$.
Moreover, $T$ is weakly pointwise dual ergodic and exhibits Darling-Kac asymptotics,
$$   \frac{S_n(f)}{a(n)}
     \overset{\mathfrak d}{\underset{n\to\infty}\lra}\
     m_T(f)\,Y_\g\ \;\;\;\; \text{for}\ f\in L_+^1(m_T),
$$
with return sequence satisfying
$  a(n) \sim c/n^{\gamma}$ for some $c>0$.
\end{propn}

\pf.
Let $\Om:=[\xi,\eta_r]$, and consider the induced map $S=T_\Om$ and the
corresponding return time function $\v=\v_\Om$. We are going to show that
the induced system naturally comes as a pcw increasing map $(\Om,S,\a)$,
which together with $\v$ satisfies the assumptions of Proposition \ref{Prop_ganzamschluss}.
Therefore Proposition \ref{T_SLT_ABU} applies, which via (\ref{eq_duality1}) entails the
DK-limit. Weak pointwise dual ergodicity implicit in the application of these
propositions.

Note first that $B_* = \cup_{n \ge 1} (\tau_{n+1},\tau_n)$, where
$\tau_1:=\xi$ and $\tau_{n+1} := w(\tau_n)$
with $w:=(T\vert_{B_*})^{-1}$ denoting the inverse of the leftmost branch of $T$.
As a consequence of (\ref{eq_indiff}) we have (Corollary on p. 82 of [T1])
\begin{equation}\label{eq_kuuen}
 q_n := \tau_n - \eta_l \sim (\kappa  n / \gamma)^{\gamma}\ \ \ \text{ as }\ n \to \infty.
\end{equation}

Fix any $B \in \b \setminus \{B_*\}$, and let $B(k):=B\cap [\v=k]$, $k\ge1$,
which defines the cylinders of $S$ inside $B$.
The induced map $S$ is trivially pcw increasing and satisfies (U).
It also satisfies (A), which is checked by the same
argument as in [T1] or [Z1], which we do not reproduce here.

Now $SB(1)=TB(1)\cap\Om \supset (\xi_l,z_T-\xi_l)$. For $k\ge2$,
$TB(k)=(\tau_{k-1},\tau_{k-2})$, and hence
$SB(k)=T^kB(k)=T(\tau_{1},\tau_0)\supset (\xi,z_T-\xi)$.
Therefore $S$ has the common image property.

Enumerating $\b \setminus \{B_*\} = \{B_0,B_1, \cdots \}$, we get
$$[\v > n] \cap B_j = v_{B_j} \left((\eta_l, \tau_{n})\right) =: (x_j,x_j+y_{j,n})$$
for suitable $x_j,y_{j,n}$,
where $v_{B_j}:=(T\vert_{B_j})^{-1}$.
This collection of intervals now
satisfies the assumptions of Lemma \ref{Lem_ManyBits} with
$s_j := v_{B_j}' (\beta_j)$ where $\beta_j$ is the left
endpoint of $B_j$. Uniformity of $y_{j,n} / q_n \to s_j$ in $j$ is a
consequence of the distortion control for (the first iterate of) $T$
provided by condition (A), which also implies $\sum s_j < \infty$.
To see that $s_j h(x_j^+) >0$ for some $j$, use the common image
property of $S$ and Lemma \ref{Lem_LowBdH} to obtain some $z_S>0$
such that $\inf_{(\xi,\xi+z_S)} h >0$. Now choose some $j$
for which $\beta_j \in (\xi,\xi+z_S)$.
\qed\\

\begin{exple}\label{expleintervalmaps}
Fix $\g \in (0,1)$ and define $F(x):=x(1+x^{1/\g})/(1-x^{1/\g})$,
$x\in X:=(0,1)$.
Let $T:X \to X$ be of the form $Tx=F(x)-F(\xi_n)$
for $x \in (\xi_n,\xi_{n+1})=:A_n$, where
$0=\xi_0 < \xi_1 < \ldots < \xi_n \nearrow 1$ are such
that $T\xi_n^- \geq T\xi_{n+1}^-$.
(It is easily seen that when $\xi_0,\ldots, \xi_n$
with these properties have been chosen, there
is a nondegenerate interval $J_n$ of admissible choices
of $\xi_{n+1}$.)
Then $T$ satisfies the assumptions of the preceding proposition,
but the finite image (F) is only fulfilled
in the exceptional cases when $T\xi_n =1$ for $n\geq n_0$
(i.e. when $\xi_{n+1}$ is the right endpoint of $J_n$ for $n\geq n_0$).
Therefore $T$ does not in general belong to the AFN maps
of [Z2].
\end{exple}

\

\end{document}